# Rust Implementation of Finite Element Exterior Calculus on Coordinate-Free Simplicial Complexes

**Abstract**

This thesis presents the development of a novel finite element library in Rust based on the principles of Finite Element Exterior Calculus (FEEC). The library solves partial differential equations formulated using differential forms on abstract, coordinate-free simplicial complexes in arbitrary dimensions, employing an intrinsic Riemannian metric derived from edge lengths via Regge Calculus. We focus on solving elliptic Hodge-Laplace eigenvalue and source problems on the $n$D de Rham complex. We restrict ourselves to first-order Whitney basis functions. The implementation is partially verified through convergence studies.

*Luis Wirth*

luwirth@ethz.ch

ethz.lwirth.com

1st May 2025

# Introduction

**Partial Differential Equations** (PDEs) are the mathematical language used to model continuum problems across physics and engineering. From heat transfer and fluid dynamics to electromagnetism and quantum mechanics, PDEs describe the fundamental systems that govern our universe. Solving PDEs efficiently and accurately is therefore central to modern computational science. [1]

The **Finite Element Method** (FEM) is one of the most important methods employed to numerically solve PDEs, particularly on unstructured meshes, drawing inspiration from **functional analysis** [1]. While scalar-valued FEM, using Lagrangian finite element spaces, is well-established, extending FEM to vector-valued problems traditionally involved intricate constructions. [2]

> Without referring to differential geometry, several authors had devised **vector valued finite elements** that can be regarded as special cases of discrete differential forms. Their constructions are formidably intricate and require much technical effort. A substantial simplification can be achieved: One should exploit the facilities provided by differential geometry for a completely coordinate-free treatment of discrete differential forms. Once we have shed the cumbersome vector calculus, everything can be constructed and studied with unmatched generality and elegance. This can be done for simplicial meshes in arbitrary dimension.
>
> — Hiptmair [2]

**Finite Element Exterior Calculus** (FEEC) developed by Douglas Arnold, Richard Falk and Ragnar Winther [3], [4] provides exactly this unified mathematical framework. By employing the language of **differential geometry** [5] and **algebraic topology** [6], FEEC extends FEM to handle problems involving **differential forms of arbitrary rank**. This approach offers robust discretizations that **preserve key topological and geometric structures** inherent in the underlying PDEs, ensuring stability, accuracy, and convergence. FEEC is now the standard framework for analyzing and constructing conforming finite element spaces for differential forms in **arbitrary dimensions** and on **domains with non-trivial topology**. [4]

A key strength of FEEC lies in its ability to naturally handle arbitrary domain topologies. This relies on fundamental connections between the algebraic topology of the simplicial complex discretizing the domain (**simplicial homology**) and the structure of differential forms on the continuous domain (**de Rham cohomology**) [6]. The de Rham theorem establishes an isomorphism, ensuring that the discrete formulation accurately captures topological features, such as **holes**, which influence the existence and uniqueness of PDE solutions. [4]

This theoretical elegance, however, is in stark contrast to many existing FEM software implementations, which are usually hard-coded to 3 dimensions and rely on vector proxies instead of actual differential forms and **exterior algebra**. [1] Furthermore, almost all implementations make use of global coordinates on the manifolds, therefore relying on embeddings instead of the intrinsic geometry of the manifold.

This thesis presents a different approach, fully embracing the **coordinate-free** perspective inherent in differential geometry for unparalleled generality. Since FEEC is formulated using differential geometry, PDE domains can be treated as **abstract Riemannian manifolds** [5]. We develop a novel finite element library in the Rust programming language [7] that operates on such abstract **simplicial complexes** [6], [8] in arbitrary dimensions, avoiding any reliance on coordinate embeddings. The geometry is defined purely intrinsically via a **Riemannian metric** [5] derived from **edge lengths** using **Regge Calculus** [9].

We restrict this work to first-order methods, employing the piecewise linear **Whitney basis** [10], [3]. This requires a **piecewise-flat approximation of the manifold geometry**, to constitute an **admissible geometric variational crime** [11].

The prototypical second-order elliptic operator in FEEC is the **Hodge-Laplace operator** [4], a generalization of the standard Laplace-Beltrami operator. Central to its analysis is **Hodge theory** [5], which provides the crucial link between the kernel of this elliptic operator (the **space of harmonic forms**) and the topology of the manifold via cohomology. Our implementation is guided by the goal of solving the Hodge-Laplace eigenvalue and source problems on the $n$D **de Rham complex** [3]. For both problems, we utilize a **mixed weak formulation** [3], [4].

The choice of **Rust** [7] stems from its strong guarantees in memory safety, performance comparable to C/C++, and modern language features suitable for complex scientific software. Its ownership model, borrow checker, and type system act like a proof system, preventing entire classes of bugs like memory errors and data races, which is crucial for reliable high-performance computing.

# Table of Contents

# Chapter 1

# Software Design & Implementation Choices

In this chapter we briefly discuss some general software engineering decisions for our library, **formoniq**.

## 1.1 Rust

We have chosen Rust [12] as the main programming language for the implementation of our Finite Element library. There are various reasons for this choice, some of which we briefly outline here.

### 1.1.1 Memory Safety + Performance

Rust is a strongly-typed, modern systems programming language that combines performance on par with C/C++ with strong memory safety guarantees. Unlike traditional memory-safe languages that rely on garbage collection, Rust achieves memory safety through a unique ownership and borrowing model inspired by formal software verification and static analysis techniques. This ensures safety at compile-time without compromising runtime performance, [12].

The Rust compiler acts as a proof checker, requiring the programmer to provide sufficient evidence for the safety of their code. This is accomplished by extending Rust's strong type system with an ownership and borrowing model that enforces clear semantics regarding memory responsibility [13].

This system completely eliminates entire classes of memory-related bugs (dangling pointers, use-after-free, data races), making software significantly more reliable. This is not limited to simple single-threaded serial computations, but extends to concurrent, parallel and distributed computing, ensuring that data races can never occur. This "fearless concurrency" feature allows developers to be fully confident that any parallel code written in Rust that compiles will behave as expected [7].

### 1.1.2 Expressiveness and Abstraction

Rust is a highly expressive language enabling powerful abstractions without sacrificing performance, often referred to as zero-cost abstractions [12]. This allows for a direct realization of concepts, making Rust well-suited for capturing precise mathematical structures naturally. Key features include:

- **Traits and Generics**: Facilitate powerful polymorphism and code reuse without the rigid hierarchies of traditional object-oriented inheritance. Compile-time monomorphization ensures zero-cost abstraction. Traits mitigate the notorious template-related compiler errors of C++, as trait bounds explicitly state required behaviors within function and struct signatures.
- **Enums, Option, and Result**: Algebraic data types, particularly enum as a tagged union, enable type-safe state representation and a simple yet powerful form of polymorphism. Option<T> eliminates null-pointer issues by enforcing explicit handling of absence, and Result<T, E> provides structured error handling without exceptional control flow.
- **Expression-Based Language and Pattern Matching**: Most constructs are expressions returning values. Combined with powerful pattern matching for destructuring composite types, this leads to concise, readable, functional-style composable code.
- **Functional Programming and Iterators**: Rust embraces higher-order functions, closures, and efficient, lazy iterators with standard methods like map, filter, and fold, promoting declarative coding.

Together, these features allow developers to write concise, maintainable, and high-performance software, combining modern paradigms with low-level control.

### 1.1.3  Infrastructure and Tooling

Beyond language features, Rust's exceptional, official tooling ecosystem streamlines development [7]. This consistency contrasts favorably with the fragmented `C++` ecosystem.

- **Cargo** is Rust's official package manager and build system, which is one of the most impressive pieces of tooling. It eliminates the need for traditional build tools like Makefiles and CMake, which are often complex and difficult to maintain—not to mention the dozens of other build systems for `C++`. Cargo simplifies dependency management through its seamless integration with `crates.io`, Rust's central package repository. Developers can effortlessly include third-party libraries by specifying them in the `Cargo.toml` file (which is versioned by git), with Cargo automatically handling downloading, compiling, and dependency resolution while enforcing semantic versioning. Publishing a crate is also straightforward.
- **Clippy** is Rust's official linter, offering valuable suggestions for improving code quality, adhering to best practices, and catching common mistakes. Our codebase does not produce a single warning or lint, passing all default checks for code quality.
- **Rustdoc** is Rust's official documentation tool, allowing developers to write inline documentation using Markdown, seamlessly integrated with code comments. This documentation can be compiled into a browsable website with `cargo doc` and is automatically published to `docs.rs` when a crate is uploaded to `crates.io`.
- **Rusttest** is the testing functionality built into Cargo for running unit and integration tests. Unit tests can be placed alongside the normal source code with a simple `#[test]` attribute, and the `cargo test` command will execute all test functions, verifying correctness. This command also ensures that all code snippets in the documentation are compiled and checked for runtime errors, keeping documentation up-to-date without requiring external test frameworks like Google Test.
- **Rustfmt** standardizes code formatting, eliminating debates about code style and ensuring consistency across projects. Our codebase fully adheres to Rustfmt's formatting guidelines. For conciseness however we will be breaking the formatting style when putting code inline into this document.

This comprehensive tooling ensures a smooth, efficient, and reliable development experience.

### 1.1.4  Drawbacks

Despite its strengths, we encountered some challenges:

- **Ecosystem Maturity**: Being relatively young (1.0 release in 2015 [12]), Rust's library ecosystem is still evolving compared to `C++`. Notably, the lack of sophisticated native sparse linear algebra solvers forced reliance on `C/C++` libraries for this project.
- **Learning Curve**: Rust's unique concepts (ownership, borrowing, lifetimes) and syntax present a steeper learning curve compared to some other languages [12].
- **Over-Engineering Risk**: The powerful type system can sometimes tempt over-engineering.

## 1.2  General software architecture

Our primary goal is to model mathematical concepts faithfully, ensuring both mathematical accuracy and code clarity. This aims to make the code mathematically expressive and self-documenting for those familiar with the underlying theory. While prioritizing mathematical rigor, we also recognize the importance of good API design and performance. Data structures are designed with efficiency and memory economy in mind, leveraging Rust's capabilities for performance-critical computing.

### 1.2.1  Modularity

FEM libraries are typically large software projects with distinct components. To promote reusability (e.g., using the mesh component in other applications) and maintainability, our library is split into multiple crates built upon each other, managed within a Cargo workspace [7].

The core crates are:
- `common`: Shared utilities and basic types.
- `manifold`: Topological and geometrical mesh data structures.

- `exterior`: Exterior algebra data structures.
- `ddf`: Discrete differential forms consisting of Cochains and Whitney forms.
- `formoniq`: The main library assembling components and providing FEEC solvers.

The chapters of the thesis parallel the structure and outline of these libraries.

## 1.3 External libraries

We briefly discuss the major external libraries used:

### 1.3.1 Nalgebra (linear algebra)

Numerical algorithms heavily rely on linear algebra libraries. Rust's `nalgebra` [14] provides an equivalent to C++'s Eigen, using generics effectively for static and dynamic matrix dimensions. It forms the foundation for nearly all numerical value manipulation in our library. We also utilize its sparse matrix capabilities via `nalgebra-sparse`.

### 1.3.2 PETSc & SLEPc (sparse solvers)

Due to the aforementioned immaturity in Rust's native sparse solver ecosystem, we utilize PETSc [15], [16], a comprehensive C/C++ library suite, for solving large sparse linear systems (specifically using its direct solvers). For the associated eigenvalue problems, we use SLEPc [17], which builds upon PETSc.

To avoid having PETSc and SLEPc as dependencies, since they are very big, we decoupled them by not using them through Rust bindings or FFI. Instead we have a separate very small PETSc/SLEPc solver program that interfaces with our Rust program through file IO. This keeps the build step of formoniq simple and nice.

### 1.3.3 Itertools (combinatorics)

Itertools [18] is a utility crate extending Rust's standard iterators with additional adaptors and methods. We use it primarily for combinatoric algorithms (permutations, combinations) essential for mesh topology and exterior algebra operations.

### 1.3.4 IndexMap (ordered HashSet)

For our skeleton data structure that contains simplices, we need a bidirectional map, between simplices and indices. For this we use `indexmap::IndexSet`, which provides set semantics while maintaining insertion order and allowing efficient lookup by index. This is provided by the `indexmap` crate [19]. This crate provides map and set data structures that maintain insertion order, which is crucial for ensuring consistent global numbering and reproducible results in our mesh data structures.

### 1.3.5 Rayon (parallelism)

To leverage multi-core processors and accelerate computationally intensive tasks, we utilize the `rayon` crate [20]. Rayon provides easy-to-use data parallelism capabilities for Rust, allowing iterators to be processed in parallel with minimal code changes. In our library, `rayon` is employed to parallelize the assembly loop, where element matrices computed independently for each cell are summed into the global sparse matrices, significantly speeding up this process on multi-core machines.

### 1.3.6 mshio (gmsh imports)

To import mesh data generated by the popular mesh generator Gmsh [21], we use the `mshio` Rust crate [22], which handles parsing of the `.msh` file format. We will use it once in our mesh implementation.

# Chapter 2

# Topology & Geometry of Simplicial Riemannian Manifolds

In this chapter, we develop various data structures and algorithms to represent and work with our Finite Element mesh. This mesh stores the topological and geometrical properties of our arbitrary-dimensional discrete PDE domain. A simplicial complex [6] will represent the topology (incidence and adjacency) of the mesh and serve as the container for all mesh entities, which are simplices. It also provides unique identification through global numbering and iteration over these entities. For the geometry, edge lengths are stored to compute the piecewise-flat (over the cells) Riemannian metric [5], using methods from Regge Calculus [9]. We also support the optional storage of global vertex coordinates if an embedding is known.

## 2.1 Coordinate Simplices

Finite Element Methods benefit from their ability to work on unstructured meshes [1]. Instead of subdividing a domain into a regular grid, FEM operates on potentially highly non-uniform meshes. The simplest mesh type suitable for such non-uniform domains are simplicial meshes. In 2D, these are the familiar triangular meshes from computer graphics, while 3D simplicial meshes are composed of tetrahedra. These building blocks must be generalized for our arbitrary-dimensional implementation.

We begin the discussion of mesh representation with a coordinate-based object that relies on an embedding in an ambient space. Later, we will discard the coordinates and rely solely on intrinsic geometry. However, for didactic purposes, starting with coordinates is helpful.

The generalization of 2D triangles and 3D tetrahedra to $n$ dimensions is called an $n$-simplex [6]. There is a type of simplex for every dimension; the first four kinds are:

- A 0-simplex is a point.
- A 1-simplex is a line segment.
- A 2-simplex is a triangle.
- A 3-simplex is a tetrahedron.

The underlying idea is that an $n$-simplex is the convex hull of $n + 1$ affinely independent points, forming the simplest possible $n$-dimensional polytope. An $n$-simplex $\sigma$ is defined by $n + 1$ vertices $\boldsymbol{v}_0, ..., \boldsymbol{v}_n \in \mathbb{R}^N$ in a possibly higher-dimensional space $\mathbb{R}^N$ (where $N \geq n$). The simplex itself is the region bounded by the convex hull of these vertices:

$$\sigma = \text{convex} \{\boldsymbol{v}_0, ..., \boldsymbol{v}_n\} = \left\{ \sum_{i=0}^{n} \lambda^i \boldsymbol{v}_i \ \middle| \ \lambda^i \geq 0, \quad \sum_{i=0}^{n} \lambda^i = 1 \right\} \tag{1}$$

We call such an object a **coordinate simplex** because it depends on the global coordinates of its vertices and resides within a potentially higher-dimensional ambient space $\mathbb{R}^N$, thus relying on an embedding. This object is uniquely determined by the vertex coordinates $\boldsymbol{v}_i$. This inspires a straightforward computational representation using a struct that stores the coordinates of each vertex as columns in a matrix:

```rust
pub struct SimplexCoords {
  pub vertices: na::DMatrix<f64>,
```

```
}
impl SimplexCoords {
  pub fn nvertices(&self) -> usize { self.vertices.ncols() }
  pub fn coord(&self, ivertex: usize) -> CoordRef { self.vertices.column(ivertex) }
}
```

We implement methods to retrieve both the intrinsic dimension $n$ (one less than the number of vertices) and the ambient dimension $N$ (the dimension of the coordinate vectors). A special and particularly simple case occurs when the intrinsic and ambient dimensions coincide, $n = N$.

```
pub type Dim = usize;
pub fn dim_intrinsic(&self) -> Dim { self.nvertices() - 1 }
pub fn dim_ambient(&self) -> Dim { self.vertices.nrows() }
pub fn is_same_dim(&self) -> bool { self.dim_intrinsic() == self.dim_ambient() }
```

### 2.1.1 Barycentric Coordinates

The coefficients $\boldsymbol{\lambda} = \left[\lambda^i\right]_{i=0}^n$ in Equation 1 are called **barycentric coordinates** [5], [1]. They appear in the convex combination $\sum_{i=0}^n \lambda^i \boldsymbol{v}_i$ as weights $\lambda^i \in [0, 1]$ (for points within the simplex) that sum to one $\sum_{i=0}^n \lambda^i = 1$. They constitute an intrinsic local coordinate representation with respect to the simplex $\sigma$, independent of the embedding in $\mathbb{R}^N$.

The coordinate transformation $\varphi : \boldsymbol{\lambda} \mapsto \boldsymbol{x}$ from intrinsic barycentric coordinates $\boldsymbol{\lambda}$ to ambient Cartesian coordinates $\boldsymbol{x}$ is given by:

$$\varphi : \boldsymbol{\lambda} \mapsto \boldsymbol{x} = \sum_{i=0}^n \lambda^i \boldsymbol{v}_i \qquad (2)$$

This transformation can be implemented as:

```
pub fn bary2global(&self, bary: CoordRef) -> Coord {
  self
    .vertices
    .coord_iter()
    .zip(bary.iter())
    .map(|(vi, baryi)| baryi * vi)
    .sum()
}
```

The barycentric coordinate representation extends beyond the simplex boundaries to the entire affine subspace spanned by the vertices. The condition $\sum_{i=0}^n \lambda^i = 1$ must still hold, but only points $\boldsymbol{x} \in \sigma$ strictly inside the simplex have all $\lambda^i \in [0, 1]$. Outside the simplex, some $\lambda^i$ will be greater than one or negative.

```
pub fn is_bary_inside(bary: CoordRef) -> bool {
  assert_relative_eq!(bary.sum(), 1.0);
  bary.iter().all(|&b| (0.0..=1.0).contains(&b))
}
```

The barycenter $\boldsymbol{m} = \frac{1}{n+1} \sum_{i=0}^n \boldsymbol{v}_i$ always has the special barycentric coordinate $\boldsymbol{\lambda} = \left[\frac{1}{n+1}\right]^{n+1}$ and gives the barycentric coordinates their name.

```
pub fn barycenter(&self) -> Coord {
  let mut barycenter = na::DVector::zeros(self.dim_ambient());
  self.vertices.column_iter().for_each(|v| barycenter += v);
  barycenter /= self.nvertices() as f64;
  barycenter
}
```

This coordinate representation treats all vertices symmetrically, assigning a weight to each. Consequently, with $n + 1$ coordinates $\lambda^0, ..., \lambda^n$ for an $n$-dimensional affine subspace subject to the constraint $\sum \lambda^i = 1$, there is redundancy. This representation is not a minimal coordinate system.

To obtain a proper coordinate system, we can single out one vertex, say $\boldsymbol{v}_0$, as the **base vertex**.

```
pub fn base_vertex(&self) -> CoordRef { self.coord(0) }
```

We can then omit the redundant coordinate $\lambda^0 = 1 - \sum_{i=1}^{n} \lambda^i$ associated with $\boldsymbol{v}_0$. The remaining **reduced barycentric coordinates** also called **local coordinates** $t = \left[\lambda^i\right]_{i=1}^{n}$ form a proper coordinate system for the $n$-dimensional affine subspace. In this system, the local coordinates $\lambda^1, ..., \lambda^n$ are unconstrained, providing a unique representation for every point in the affine subspace via a bijection with $\mathbb{R}^n$.

```rust
pub fn bary2local(bary: CoordRef) -> Coord {
  bary.view_range(1.., ..).into()
}
pub fn local2bary(local: CoordRef) -> Coord {
  let bary0 = 1.0 - local.sum();
  local.insert_row(0, bary0)
}
```

## 2.1.2 Spanning Vectors

Consider the edge vectors emanating from the base vertex: $\boldsymbol{e}_i = \boldsymbol{v}_i - \boldsymbol{v}_0 \in \mathbb{R}^N$ for $i = 1, ..., n$. These are the **spanning vectors**. We can collect them as columns into a matrix $\mathbf{E} \in \mathbb{R}^{N \times n}$:

$$\mathbf{E} = \begin{bmatrix} | & & | \\ \boldsymbol{e}_1 & \cdots & \boldsymbol{e}_n \\ | & & | \end{bmatrix} \tag{3}$$

This matrix can be computed as follows:

```rust
pub fn spanning_vectors(&self) -> na::DMatrix<f64> {
  let mut mat = na::DMatrix::zeros(self.dim_ambient(), self.dim_intrinsic());
  let v0 = self.base_vertex();
  // Skip base vertex (index 0)
  for (i, vi) in self.vertices.column_iter().skip(1).enumerate() {
    let v0i = vi - v0;
    mat.set_column(i, &v0i);
  }
  mat
}
```

This gives us an explicit basis of the affine space located at $\boldsymbol{v}_0$ and spanned by $\boldsymbol{e}_1, ..., \boldsymbol{e}_n$.

We can also rewrite the coordinate transformation $\varphi : \boldsymbol{\lambda} \mapsto \boldsymbol{x}$ using $\lambda^0 = 1 - \sum_{i=1}^{n} \lambda^i$ in terms of the spanning vectors instead of the vertices:

$$\boldsymbol{x} = \sum_{i=0}^{n} \lambda^i \boldsymbol{v}_i = \left(1 - \sum_{i=1}^{n} \lambda^i\right) \boldsymbol{v}_0 + \sum_{i=1}^{n} \lambda^i \boldsymbol{v}_i = \boldsymbol{v}_0 + \sum_{i=1}^{n} \lambda^i (\boldsymbol{v}_i - \boldsymbol{v}_0) = \boldsymbol{v}_0 + \mathbf{E}\boldsymbol{t} \tag{4}$$

This shows that the transformation $\varphi$ is actually an affine map $\varphi : \boldsymbol{t} \mapsto \boldsymbol{x}$ consisting of the linear map represented by $\mathbf{E}$ followed by a translation by $\boldsymbol{v}_0$.

$$\varphi : \boldsymbol{t} \mapsto \boldsymbol{x} = \boldsymbol{v}_0 + \mathbf{E}\boldsymbol{t} \tag{5}$$

We implement functions for this affine transformation:

```rust
pub fn linear_transform(&self) -> na::DMatrix<f64> { self.spanning_vectors() }
pub fn affine_transform(&self) -> AffineTransform {
  let translation = self.base_vertex().into_owned();
  let linear = self.linear_transform();
  AffineTransform::new(translation, linear)
}
pub fn local2global(&self, local: CoordRef) -> Coord {
  self.affine_transform().apply_forward(local)
}
```

This makes use of the `AffineTransform` struct and its forward application:

```rust
pub struct AffineTransform {
  pub translation: na::DVector<f64>,
  pub linear: na::DMatrix<f64>,
```

```rust
}
impl AffineTransform {
  pub fn new(translation: Vector, linear: Matrix) -> Self { Self { translation, linear } }
  pub fn dim_domain(&self) -> usize { self.linear.ncols() }
  pub fn dim_image(&self) -> usize { self.linear.nrows() }
  pub fn apply_forward(&self, coord: na::DVectorView<f64>) -> na::DVector<f64> {
    &self.linear * coord + &self.translation
  }
}
```

The reverse transformation $\psi : \boldsymbol{x} \mapsto \boldsymbol{t}$ from global $\boldsymbol{x}$ to local $\boldsymbol{t}$ coordinates is more complex due to the potentially higher-dimensional ambient space $\mathbb{R}^N$ with $N \geq n$. The global coordinate point might not lie exactly in the affine subspace due to floating-point inaccuracies. This makes the linear system $\mathbf{E}\boldsymbol{t} = \boldsymbol{x} - \boldsymbol{v}_0$ for the reverse transformation potentially underdetermined. We use the Moore-Penrose pseudo-inverse $\mathbf{E}^\dagger$ [23], typically computed via Singular Value Decomposition (SVD), to find the unique least-squares solution of smallest norm:

$$\psi : \boldsymbol{x} \mapsto \boldsymbol{t} = \mathbf{E}^\dagger (\boldsymbol{x} - \boldsymbol{v}_0) \tag{6}$$

```rust
impl SimplexCoords {
  pub fn global2local(&self, global: CoordRef) -> Coord {
    self.affine_transform().apply_backward(global)
  }
}

impl AffineTransform {
  pub fn apply_backward(&self, coord: Coord) -> Vector {
    if self.dim_domain() == 0 { return Coord::zeros(0); }
    self
      .linear
      .clone()
      .svd(true, true)
      .solve(&(coord - &self.translation), 1e-12)
      .expect("SVD solve failed")
  }
  pub fn pseudo_inverse(&self) -> Self {
    if self.dim_domain() == 0 {
      return Self::new(Coord::zeros(0), Matrix::zeros(0, self.dim_image()));
    }
    let linear = self.linear.clone().pseudo_inverse(1e-12).unwrap();
    let translation = &linear * &self.translation;
    Self { translation, linear }
  }
}
```

The derivatives of the affine transformation $\varphi : \boldsymbol{t} \mapsto \boldsymbol{x}$ reveal that the spanning vectors $\boldsymbol{e}_i$ form a natural basis for the tangent space $T_p \sigma$ at any point $p$ within the simplex $\sigma$ [5]. The Jacobian of the affine map is precisely $\mathbf{E}$.

$$\frac{\partial}{\partial \lambda^i} = \frac{\partial \boldsymbol{x}}{\partial \lambda^i} = \boldsymbol{e}_i \qquad \frac{\partial \boldsymbol{x}}{\partial \boldsymbol{t}} = \mathbf{E} \tag{7}$$

In differential geometry terms, the linear map represented by $\mathbf{E}$ acts as the **pushforward**.

$$\varphi_* : \boldsymbol{u} \mapsto \boldsymbol{w} = \mathbf{E}\boldsymbol{u} \tag{8}$$

It transforms intrinsic tangent vectors $\boldsymbol{u} = u^i \frac{\partial}{\partial \lambda^i}$ to ambient tangent vectors $\boldsymbol{w} = w^i \boldsymbol{e}_i$.

```rust
/// Local2Global Tangentvector
pub fn pushforward_vector(&self, local: TangentVectorRef) -> TangentVector {
  self.linear_transform() * local
}
```

Conversely, the **pullback** $\varphi^*$ operation takes covectors defined in the ambient space coordinates and expresses them in the local coordinate system. If $\boldsymbol{\omega} \in \mathbb{R}^{1 \times N}$ is a covector in the ambient space, its pullback $\boldsymbol{\eta} = \varphi^* \boldsymbol{\omega}$ acts on local tangent vectors $\boldsymbol{u}$ such that $\boldsymbol{\eta}(\boldsymbol{u}) = \boldsymbol{\omega}(\varphi_* \boldsymbol{u}) = \boldsymbol{\omega}(\mathbf{E}\boldsymbol{u})$. The transformation rule is right-multiplication if covectors are row vectors.

$$\varphi^* : \boldsymbol{\omega} \mapsto \boldsymbol{\eta} = \boldsymbol{\omega}\mathbf{E} \tag{9}$$

```rust
/// Global2Local Cotangentvector
pub fn pullback_covector(&self, global: CoTangentVectorRef) -> CoTangentVector {
  global * self.linear_transform()
}
```

Separately, we can consider the differentials of the barycentric coordinate functions $\lambda^i$ as functions of the global coordinates $\boldsymbol{x}$. These differentials, $\mathrm{d}\lambda^i$, form a basis for the cotangent space $T_p^*\sigma$. Their components relative to the ambient basis $\mathrm{d}x^i$ are found using the differential of the inverse map $\psi : \boldsymbol{x} \to \boldsymbol{t}$, which involves the pseudo-inverse $\mathbf{E}^\dagger$. Specifically, the rows of $\mathbf{E}^\dagger$ give the components of $\mathrm{d}\lambda^1, ..., \mathrm{d}\lambda^n$. The differential $\mathrm{d}\lambda^0$ is determined by the constraint $\sum_i \mathrm{d}\lambda^i = 0$.

$$\frac{\partial \boldsymbol{t}}{\partial \boldsymbol{x}} = \mathbf{E}^\dagger \qquad \mathrm{d}\lambda^i = \frac{\partial \lambda^i}{\partial \boldsymbol{x}} = \varepsilon^i = \left(\mathbf{E}^\dagger\right)_{i,:} \qquad \mathrm{d}\lambda^0 = -\sum_{i=1}^n \mathrm{d}\lambda^i \tag{10}$$

These $\varepsilon^i = \mathrm{d}\lambda^i$ form the basis dual to the tangent basis $\boldsymbol{e}_1, ..., \boldsymbol{e}_n$ [5].

$$\mathrm{d}\lambda^i\left(\frac{\partial}{\partial \lambda^j}\right) = \delta_j^i \tag{11}$$

```rust
/// Total differential of barycentric coordinate functions in the rows(!) of
/// a matrix.
pub fn difbarys(&self) -> Matrix {
  let difs = self.inv_linear_transform();
  let mut difs = difs.insert_row(0, 0.0);
  difs.set_row(0, &-difs.row_sum()); // lambda^0 = -sum(dif lambda^i)
  difs
}
```

The spanning vectors as a basis of the tangent space, can be used to derive the **Riemannian metric tensor**. This metric on the simplex is induced by the ambient Euclidean metric via the affine embedding.

$$\mathbf{G} = \mathbf{E}^\top\mathbf{E} \qquad \mathbf{G}_{ij} = e_i \cdot e_j \tag{12}$$

Our implementation has a struct for working with Gramians, which will discuss in more detail later on.

```rust
pub fn metric_tensor(&self) -> Gramian {
  Gramian::from_euclidean_vectors(self.spanning_vectors())
}
```

Since $\mathbf{E}$ is constant across the simplex, this induced metric $\mathbf{G}$ is also **constant** everywhere within the simplex. This constancy implies that the simplex is **intrinsically flat** (zero Riemannian curvature).

The spanning vectors also define a parallelepiped. The volume of the $n$-simplex $\sigma$ is $\frac{1}{n!}$ times the $n$-dimensional volume of this parallelepiped. The signed volume is computed using the determinant of the spanning vectors if $n = N$.

$$\mathrm{vol}(\sigma) = \frac{1}{n!}\det(\mathbf{E}) \tag{13}$$

For higher-dimensional ambient spaces we can use the Gram determinant of the metric tensor. [5].

$$\mathrm{vol}(\sigma) = \frac{1}{n!}\sqrt{\mathbf{G}} = \frac{1}{n!}\sqrt{\det(\mathbf{E}^\top\mathbf{E})} \tag{14}$$

```rust
impl SimplexCoords {
  pub fn det(&self) -> f64 {
    let det = if self.is_same_dim() {
      self.spanning_vectors().determinant()
    } else {
      self.metric_tensor().det_sqrt()
    };
    refsimp_vol(self.dim_intrinsic()) * det
  }
```

```rust
  pub fn vol(&self) -> f64 { self.det().abs() }
  pub fn is_degenerate(&self) -> bool { self.vol() <= 1e-12 }
}
pub fn refsimp_vol(dim: Dim) -> f64 { factorialf(dim).recip() }
```

The sign of the signed volume gives the global orientation of the coordinate simplex relative to the ambient space.

```rust
pub fn orientation(&self) -> Sign {
  Sign::from_f64(self.det()).unwrap()
}

#[derive(Debug, Clone, Copy, Default, PartialEq, Eq, Hash)]
pub enum Sign {
  #[default]
  Pos = 1,
  Neg = -1,
}
impl Sign {
  pub fn from_f64(f: f64) -> Option<Self> {
    if f == 0.0 {
      return None;
    }
    Some(Self::from_bool(f > 0.0))
  }
}
```

Swapping two vertices in the definition of the simplex negates the determinant and thus flips the orientation. Every non-degenerate simplex has exactly two orientations, positive and negative [6].

| Simplex Orientation | | | |
|---|---|---|---|
| $n$ | Simplex | Positive | Negative |
| 1 | Line Segment | left-to-right | right-to-left |
| 2 | Triangle | counterclockwise | clockwise |
| 3 | Tetrahedron | right-handed | left-handed |

### 2.1.3 Reference Simplex

A particularly important simplex is the **reference simplex** [1], which serves as a canonical domain for defining basis functions in FEM. For a given dimension $n$, the $n$-dimensional reference simplex $\hat{\sigma}^n$ is defined in $\mathbb{R}^n$ (so $N = n$) using its local coordinates as its global Cartesian coordinates:

$$\hat{\sigma}^n = \left\{ (\lambda_1, ..., \lambda_n) \in \mathbb{R}^n \,\middle|\, \lambda_i \geq 0, \quad \sum_{i=1}^{n} \lambda_i \leq 1 \right\} \tag{15}$$

Its vertices are the origin $\boldsymbol{v}_0 = \boldsymbol{0}$ and the standard basis vectors $\boldsymbol{v}_i = \hat{\boldsymbol{e}}_i$ for $i = 1...n$. The spanning vectors are simply the standard basis vectors $\boldsymbol{e}_i = \hat{\boldsymbol{e}}_i$, so $\mathbf{E} = \mathbf{I}_n$. The metric tensor is the identity matrix $\mathbf{G} = \mathbf{I}_n$, representing the standard Euclidean inner product. Its volume is $(n!)^{-1}$.

When looking at an arbitrary "real" $n$-simplex $\tau$, the local to global map $\varphi_\tau : \boldsymbol{t} \mapsto \boldsymbol{x}$ can be seen as a parametrization $\varphi_\tau : \hat{\sigma}^n \to \tau \subseteq \mathbb{R}^N$, where the parametrization domain is the reference $n$-simplex. The real simplex is then the image of the reference simplex $\sigma = \varphi_\sigma(\hat{\sigma}^n)$

Conversely the global to local map $\psi_\tau : \boldsymbol{x} \mapsto \boldsymbol{t}$ can be seen as a chart map $\psi_\tau : \tau \to \hat{\sigma}^n$ where the chart itself is the reference $n$-simplex.

## 2.2 Abstract Simplices

After studying coordinate simplices, the reader has hopefully developed some intuitive understanding of simplices. We will now shed the coordinates and represent simplices in a more abstract way, by considering them merely as a list of vertex indices, without any associated vertex coordinates. An $n$-simplex $\sigma$ is thus represented as an $(n + 1)$-tuple of natural numbers, which correspond to vertex indices [6].

$$\sigma = [v_0, ..., v_n] \in \mathbb{N}^{n+1} \qquad v_i \in \mathbb{N} \tag{16}$$

In Rust, we can simply represent this using the following struct:

```rust
pub type VertexIdx = usize;
pub struct Simplex {
  pub vertices: Vec<VertexIdx>,
}
impl Simplex {
  pub fn new(vertices: Vec<VertexIdx>) -> Self { Self { vertices } }
  pub fn standard(dim: Dim) -> Self { Self::new((0..dim + 1).collect()) }
  pub fn single(v: usize) -> Self { Self::new(vec![v]) }

  pub fn nvertices(&self) -> usize { self.vertices.len() }
  pub fn dim(&self) -> Dim { self.nvertices() - 1 }
}
```

The ordering of the vertices *does* matter; therefore, we are dealing with ordered tuples, not just unordered sets. This makes our simplices combinatorial objects, and these combinatorics will be the heart of our mesh data structure.

### 2.2.1 Sorted Simplices

Even though the order matters for defining a specific simplex instance, simplices that share the same set of vertices are closely related. For this reason, it is helpful to introduce a convention for the canonical representation of a simplex given a particular set of vertices.

Our canonical representative will be the tuple whose vertex indices are sorted in increasing order. We can take any simplex with arbitrarily ordered vertices and convert it to its canonical representation.

```rust
pub fn is_sorted(&self) -> bool { self.vertices.is_sorted() }
pub fn sort(&mut self) { self.vertices.sort_unstable() }
pub fn sorted(mut self) -> Self { self.sort(); self }
```

Using this canonical representation, we can easily check whether two simplices have the same vertex set, meaning they are permutations of each other.

```rust
pub fn set_eq(&self, other: &Self) -> bool {
  self.clone().sorted() == other.clone().sorted()
}
pub fn is_permutation_of(&self, other: &Self) -> bool { self.set_eq(other) }
```

### 2.2.2 Orientation

For coordinate simplices, we observed that a simplex always has two possible orientations. We computed this orientation based on the determinant of the spanning vectors, but without coordinates, this is no longer possible.

However, we can still define a notion of relative orientation. Recall that swapping two vertices in a coordinate simplex flips its orientation due to the properties of the determinant. The same behavior is present in abstract simplices, based purely on the ordering of the vertices. All permutations of a given set of vertices can be divided into two equivalence classes relative to a reference ordering: even and odd permutations. We associate simplices with even permutations with positive orientation and those with odd permutations with negative orientation. Therefore, every abstract simplex has exactly two orientations, positive and negative, depending on its vertex ordering relative to a reference. [6]

We use our canonical sorted representation as the reference ordering. We can determine the orientation of any simplex relative to this sorted permutation by counting the number of swaps necessary to sort its vertex list. An even number of swaps corresponds to a positive orientation, and an odd number corresponds to a negative orientation. For this, we implement a basic bubble sort that keeps track of the number of swaps. The computational complexity $\mathcal{O}(n^2)$ is not optimal but sufficient for the typically small number of vertices per simplex.

```rust
pub fn orientation_rel_sorted(&self) -> Sign { self.clone().sort_signed() }
pub fn sort_signed(&mut self) -> Sign { sort_signed(&mut self.vertices) }
```

```rust
/// Returns the sorted permutation of `a` and the sign of the permutation.
pub fn sort_signed<T: Ord>(a: &mut [T]) -> Sign {
  Sign::from_parity(sort_count_swaps(a))
}
/// Returns the sorted permutation of `a` and the number of swaps.
pub fn sort_count_swaps<T: Ord>(a: &mut [T]) -> usize {
  let mut nswaps = 0;

  let mut n = a.len();
  if n > 0 {
    let mut swapped = true;
    while swapped {
      swapped = false;
      for i in 1..n {
        if a[i - 1] > a[i] {
          a.swap(i - 1, i);
          swapped = true;
          nswaps += 1;
        }
      }
      n -= 1;
    }
  }
  nswaps
}
```

Two simplices composed of the same vertex set have equal orientation if and only if their respective permutations fall into the same equivalence class (both even or both odd). Using the transitivity of this equivalence relation, we can check this by comparing their orientations relative to the canonical sorted permutation.

```rust
pub fn orientation_eq(&self, other: &Self) -> bool {
  self.orientation_rel_sorted() == other.orientation_rel_sorted()
}
```

### 2.2.3 Subsets

Another important notion is that of a subsimplex or a face of a simplex [6]. In this context, we first consider the concept of subset simplices.

A simplex $\sigma$ can be considered as the subset of another simplex $\tau$ if all vertices of $\sigma$ are also vertices of $\tau$. We can check this condition directly: $\sigma \subseteq \tau \Leftrightarrow (\forall a \in \sigma \Rightarrow a \in \tau)$

```rust
pub fn is_subset_of(&self, other: &Self) -> bool {
  self.iter().all(|v| other.vertices.contains(v))
}
pub fn is_superset_of(&self, other: &Self) -> bool {
  other.is_subset_of(self)
}
```

Note that this only considers the simplices as sets, even though they are actually lists with an ordering.

Given a simplex that is a subset of a larger simplex, we can compute its vertex indices relative to:

```rust
/// Computes local vertex numbers relative to sup.
pub fn relative_to(&self, sup: &Self) -> Simplex {
  let local = self
    .iter()
    .map(|iglobal| {
      sup
        .iter()
        .position(|iother| iglobal == iother)
        .expect("Not a subset.")
    })
    .collect();
  Simplex::new(local)
}
```

### 2.2.4   Subsequences

When considering the subsimplices of an $n$-simplex there are multiple simplices (as ordered vertex lists) that have the same vertex set. However it's desirable to have a unique representative simplex for each vertex subset. A natural choice for this is the subset that preserves the ordering of the vertices. This is particularly important for the boundary operator.

For this, we consider **subsequences** of the original simplex's vertex list. A subsequence maintains the relative order of the vertices it contains. We can check if one simplex is a subsequence of another using a naive algorithm:

```rust
pub fn is_subsequence_of(&self, sup: &Self) -> bool {
  let mut sup_iter = sup.iter();
  self.iter().all(|item| sup_iter.any(|x| x == item))
}
pub fn is_supersequence_of(&self, other: &Self) -> bool {
  other.is_subsequence_of(self)
}
```

We provide a method for generating all $k$-dimensional subsimplices that are **subsequences** of an $n$-simplex. This is achieved by generating all $(k + 1)$-length selections of the original $(n + 1)$ vertices, using an implementation provided by the itertools crate [18].

```rust
pub fn subsequences(&self, sub_dim: Dim) -> impl Iterator<Item = Self> {
  itertools::Itertools::combinations(self.clone().into_iter(), sub_dim + 1).map(Self::from)
}
```

This implementation conveniently provides the subsequences in lexicographical order with respect to the vertex indices of the original simplex. If the original simplex was sorted, then the generated subsequences are also absolutely lexicographically ordered.

A standard operation is to generate all subsequence simplices of the standard simplex $[0, 1, ..., n]$. We call these the standard subsimplices.

```rust
pub fn standard_subsimps(dim_cell: Dim, dim_sub: Dim) -> impl Iterator<Item = Simplex> {
  Simplex::standard(dim_cell).subsequences(dim_sub)
}
```

We can also generate all standard subsimplices for each dimension up to $n$ in a graded fashion, which is useful for generating a standard simplicial complex.

```rust
pub fn graded_subsimps(dim_cell: Dim) -> impl Iterator<Item = impl Iterator<Item = Simplex>> {
  (0..=dim_cell).map(move |d| standard_subsimps(dim_cell, d))
}
```

Conversely, we can generate supersequences. Given a simplex and a "root" simplex (which contains both the original simplex and its potential supersequences as subsequences), we can find all subsequences of the root that have a specific `super_dim` and contain the original simplex as a subsequence.

```rust
pub fn supersequences(
  &self,
  super_dim: Dim,
  root: &Self,
) -> impl Iterator<Item = Self> + use<'_> {
  root
    .subsequences(super_dim)
    .filter(|sup| self.is_subsequence_of(sup))
}
```

The number of distinct vertex subsequences of size $k + 1$ within a set of $n + 1$ vertices is given by the binomial coefficient $\binom{n+1}{k+1}$.

```rust
pub fn nsubsimplices(dim_cell: Dim, dim_sub: Dim) -> usize {
  binomial(dim_cell + 1, dim_sub + 1)
}
```

### 2.2.5   Boundary

A special operation related to subsequence simplices is the **boundary operator** $\partial$ [6]. Applied to an $n$-simplex $\sigma = [v_0, ..., v_n]$, the boundary operator is defined as the formal sum:

$$\partial\sigma = \sum_{i=0}^{n} (-1)^i [v_0, ..., \hat{v}_i, ..., v_n] \tag{17}$$

Here, $\hat{v}_i$ indicates that vertex $v_i$ is omitted. The result is a formal sum of all $(n-1)$-dimensional subsequence simplices of $\sigma$. Each subsimplex is assigned a sign $(-1)^i$, giving the boundary an orientation consistent with the orientation of $\sigma$. When viewed as elements of the free Abelian group generated by all oriented simplices, this operator is linear.

For instance, the boundary of the triangle $\sigma = [0, 1, 2]$ is

$$\partial\sigma = (-1)^0 [1, 2] + (-1)^1 [0, 2] + (-1)^2 [0, 1] = [1, 2] - [0, 2] + [0, 1] \tag{18}$$

Rearranging terms to follow a path gives $[0, 1] + [1, 2] + [2, 0]$, which corresponds to traversing the edges of the triangle.

Our implementation generates the boundary facets using the subsequences method, which yields them in an order based on the index of the retained vertices (lexicographical combinations). This order differs from the summation index $i$ (based on the omitted vertex) in the standard definition $\sum_i (-1)^i [...]$. The code calculates the appropriate sign for each generated subsequence to ensure consistency with the standard alternating sign convention of the boundary operator.

```rust
pub fn boundary(&self) -> impl Iterator<Item = SignedSimplex> {
  let mut sign = Sign::from_parity(self.nvertices() - 1);
  self.subsequences(self.dim() - 1).map(move |simp| {
    let this_sign = sign;
    sign.flip();
    SignedSimplex::new(simp, this_sign)
  })
}
```

To represent the terms in the formal sum, we introduce a struct that pairs a simplex with an explicit sign:

```rust
#[derive(Debug, Default, Clone, PartialEq, Eq, Hash)]
pub struct SignedSimplex {
  pub simplex: Simplex,
  pub sign: Sign,
}
```

## 2.3  Simplicial Skeleton

Simplices are the building blocks of our mesh. If we construct our mesh from coordinate simplices, we obtain the typical Euclidean extrinsic description of an embedded mesh, which inherently contains all geometric information. As an embedding, the union of these coordinate simplices forms an $n$-dimensional region within an $N$-dimensional Euclidean space.

If we, in contrast, build our mesh using only abstract simplices, we lack this explicit geometric information, since abstract simplices only specify the **connectivity** via vertices indices. This information, however, fully defines the topology of our discrete $n$-manifold [6]. The connectivity (incidence and adjacency) between simplices is determined entirely by shared vertices. This makes the topology purely combinatorial.

To define the topology of a simplicial $n$-manifold at its highest dimension, we only need to store the set of $n$-simplices that constitute it. This collection defines the top-level structure of the mesh. Such a collection of $n$-simplices, typically sharing vertices from a common pool, is called an $n$-skeleton

```rust
/// A container for sorted simplices of the same dimension.
#[derive(Default, Debug, Clone)]
pub struct Skeleton {
  /// Every simplex is sorted.
  simplices: IndexSet<Simplex>,
  nvertices: usize,
}
```

A `Skeleton` fulfills several responsibilities of a mesh data structure. Primarily, it serves as a container for all $n$-simplices of a specific dimension $n$. It allows for iteration over all these mesh entities:

```rust
impl Skeleton {
  pub fn iter(&self) -> indexmap::set::Iter<'_, Simplex> {
    self.simplices.iter()
  }
}
impl IntoIterator for Skeleton {
  type Item = Simplex;
  type IntoIter = indexmap::set::IntoIter<Self::Item>;
  fn into_iter(self) -> Self::IntoIter {
    self.simplices.into_iter()
  }
}
```

Crucially, it provides unique identification for each simplex through a global numbering scheme within that dimension. This establishes a bijective mapping between an integer index KSimplexIdx and the abstract simplex itself. This functionality is achieved using the `IndexSet` data structure from the `indexmap` crate. `IndexSet` maintains insertion order, allowing efficient retrieval of a `Simplex` given its index (similar to `Vec`), while also using hashing internally to support the reverse lookup: retrieving the index corresponding to a given `Simplex` instance.

```rust
pub fn simplex_by_kidx(&self, idx: KSimplexIdx) -> &Simplex {
  self.simplices.get_index(idx).unwrap()
}
pub fn kidx_by_simplex(&self, simp: &Simplex) -> KSimplexIdx {
  self.simplices.get_index_of(simp).unwrap()
}
```

The `Skeleton` constructor enforces several guarantees about the simplices it contains:

```rust
pub fn new(simplices: Vec<Simplex>) -> Self {
  assert!(!simplices.is_empty(), "Skeleton must not be empty");
  let dim = simplices[0].dim();
  assert!(
    simplices.iter().map(|simp| simp.dim()).all(|d| d == dim),
    "Skeleton simplices must have same dimension."
  );
  assert!(
    simplices.iter().all(|simp| simp.is_sorted()),
    "Skeleton simplices must be sorted."
  );
  let nvertices = if dim == 0 {
    assert!(simplices.iter().enumerate().all(|(i, simp)| simp[0] == i));
    simplices.len()
  } else {
    simplices
      .iter()
      .map(|simp| simp.iter().max().expect("Simplex is not empty."))
      .max()
      .expect("Simplices is not empty.")
      + 1
  };

  let simplices = IndexSet::from_iter(simplices);
  Self {
    simplices,
    nvertices,
  }
}
```

First, a skeleton cannot be empty, and all contained simplices must have the same dimension. Furthermore, we enforce that only the canonical representation (sorted vertex indices) of simplices is stored. This is essential for the reverse mapping (simplex-to-index lookup) to be consistently useful; regardless of the initial ordering of a simplex's vertices, converting it to the canonical sorted form allows retrieval of its unique index within the skeleton. Lastly, a special requirement applies to 0-skeletons: the simplices must represent the vertices indexed sequentially as $[0], [1], ..., [N-1]$. The constructor also determines and stores the total number of vertices (`nvertices`) involved in the skeleton.

## 2.4 Simplicial Complex

An $n$-skeleton provides the top-level topological information of a mesh by defining its $n$-dimensional cells $\Delta_n(\mathcal{M})$. However, this structure alone is often insufficient for applications like FEEC [4]. This is because degrees of freedom (DOFs) or basis functions are not only associated with cells but also with the lower-dimensional subsimplices $\Delta_k(\mathcal{M}), k \leq n$. Therefore, we need a data structure that explicitly represents the topology of *all* relevant simplices within the mesh, including subsimplices.

Enter the **simplicial complex** [6]. A simplicial complex $K$ is a collection of simplices such that:
1. Every subsimplex of a simplex in $K$ is also in $K$.
2. The intersection of any two simplices in $K$ is either empty or a face of both.

In our implementation, we represent an $n$-dimensional simplicial complex by storing the complete set of $k$-simplices for each dimension $k$ from $0$ to $n$. Effectively, it comprises $n+1$ distinct skeletons, one for each dimension.

```
#[derive(Default, Debug, Clone)]
pub struct Complex {
  skeletons: Vec<ComplexSkeleton>,
}
impl Complex {
  pub fn dim(&self) -> Dim { self.skeletons.len() - 1 }
}

/// A skeleton inside of a complex, pairing the raw Skeleton
/// with additional topological data computed during complex construction.
#[derive(Default, Debug, Clone)]
pub struct ComplexSkeleton {
  skeleton: Skeleton,
  complex_data: SkeletonComplexData,
}

pub type SkeletonComplexData = Vec<SimplexComplexData>;

/// Complex-specific data associated with a single simplex within the complex.
#[derive(Default, Debug, Clone)]
pub struct SimplexComplexData {
  /// Stores the indices of the top-level cells (n-simplices)
  /// that contain this simplex as a face (subsequence).
  pub cocells: Vec<SimplexIdx>,
}
```

We store some minimal additional information to the skeletons, which are the parent cells (cocells) of each simplex. This allows us to find the supersimplices of each simplex, because then we have a root for searching in. This allows us to traverse the whole hierarchy of simplices efficiently.

We use the following terminology for simplices of various dimensions within the complex:
- The 0-simplices are called **vertices**.
- The 1-simplices are called **edges**.
- The 2-simplices are called **faces**.
- The 3-simplices are called **tets**.
- The $(n-1)$-simplices are called **facets**.
- The $n$-simplices are called **cells**.

The complex will serve as the main **topological data structure** passed as an argument into our FEEC algorithm routines, providing access to all mesh entities and their relationships.

While general simplicial complexes can represent complex topological spaces, potentially including non-manifold features [6], our target applications in PDE modeling typically require computational domains that are **manifolds**. A topological $n$-manifold is a space that locally resembles Euclidean $n$-space.

Our `Complex` data structure is designed to only represent such manifold domains. We ensure two properties through construction: [8]

1. **Closure:** The complex contains all faces (subsequences) of its simplices. This is guaranteed by generating the complex from a list of top-level cells and explicitly adding all their subsequences.

2. **Purity:** Every simplex of dimension $k < n$ is a face of at least one $n$-simplex (cell). This is also ensured by constructing the complex downwards from the cells.

The primary method for creating a `Complex` is by providing the skeleton of its highest-dimensional cells ($n$-simplices). The `from_cells` constructor then systematically builds the skeletons for all lower dimensions ($k = 0, ..., n - 1$) by generating all subsequence simplices of the input cells. During this process, it also stores the cocells of each simplex.

However, the input cells skeleton itself might implicitly define a non-manifold topology. A common example in 2D is when three or more triangles meet along a single edge, like pages of a book. To ensure the represented space is a manifold (potentially with boundary), we perform a check. A necessary condition for an $n$-dimensional simplicial complex to be a manifold is that every facet must be shared by at most 2 cells. If facets are shared by only one cell, this indicates they lie on the boundary of the manifold. This check is performed after building the full complex structure because the required incidence information (which cells contain each facet) is naturally computed during the construction process. It fundamentally serves as a validation step for the input cell skeleton.

```rust
impl Complex {
  pub fn from_cells(cells: Skeleton) -> Self {
    let dim = cells.dim();

    let mut skeletons = vec![ComplexSkeleton::default(); dim + 1];
    skeletons[0] = ComplexSkeleton {
      skeleton: Skeleton::new((0..cells.nvertices()).map(Simplex::single).collect()),
      complex_data: (0..cells.nvertices())
        .map(|_| SimplexComplexData::default())
        .collect(),
    };

    for (icell, cell) in cells.iter().enumerate() {
      for (
        dim_skeleton,
        ComplexSkeleton {
          skeleton,
          complex_data: mesh_data,
        },
      ) in skeletons.iter_mut().enumerate()
      {
        for sub in cell.subsequences(dim_skeleton) {
          let (sub_idx, is_new) = skeleton.insert(sub);
          let sub_data = if is_new {
            mesh_data.push(SimplexComplexData::default());
            mesh_data.last_mut().unwrap()
          } else {
            &mut mesh_data[sub_idx]
          };
          sub_data.cocells.push(SimplexIdx::new(dim, icell));
        }
      }
    }

    // Topology checks.
    if dim >= 1 {
      let facet_data = skeletons[dim - 1].complex_data();
      for SimplexComplexData { cocells } in facet_data {
        let nparents = cocells.len();
        let is_manifold = nparents == 2 || nparents == 1;
        assert!(is_manifold, "Topology must be manifold.");
      }
    }
```

```
    Self { skeletons }
  }
}
```

## 2.4.1 Simplices in the Mesh: Simplex Indices and Handles

To identify a simplex inside the mesh, we use an indexing system. If the context of a concrete dimension is given, then we only need to know the index inside the skeleton, which is just an integer

```
pub type KSimplexIdx = usize;
```

If the skeleton context is not given, then we also need to specify the dimension, for this we have a fat index, that contains both parts.

```
#[derive(Debug, Clone, Copy, PartialEq, Eq, Hash)]
pub struct SimplexIdx {
  pub dim: Dim,
  pub kidx: KSimplexIdx,
}
```

Using this we reference any simplex in the mesh. Since we are using a `indexmap::IndexSet` data structure for storing our simplices inside a skeleton, we are able to go both ways. We can not only retrieve the simplex, given the index, but we can also get the index corresponding to a given simplex.

The `Simplex` struct doesn't reference the mesh and therefore doesn't have access to any other simplices. But for doing any kind of topological computations, it is helpful to be able to reference other simplices in the mesh. For this reason we introduce a new concept, that represents a simplex inside of a mesh. We create a simplex handle, that is like a more sophisticated pointer to a simplex, that has a reference to the mesh. This allows us to interact with these simplices inside the mesh very naturally.

```
#[derive(Copy, Clone)]
pub struct SimplexHandle<'c> {
  complex: &'c Complex,
  idx: SimplexIdx,
}
```

Pretty much all functionality that was available on the raw `Simplex` struct, is also available on the handle, but here we directly provide any relevant mesh information.

For instance the `SimplexHandle::supersimps` gives us all supersimplices that are also contained in the mesh, which are exactly all the supersequence simplices. The `Simplex::supersequence` method however expects a `root: &Simplex`, which gives the context in which we are searching for supersequences. This context is directly provided to the method on the handle, since the mesh knows the cocells of each simplex, which is here chosen as the root, so we don't need to provide this argument ourselves.

```
  pub fn supersimps(&self, dim_super: Dim) -> Vec<SimplexHandle> {
    self
      .cocells()
      .flat_map(|parent| {
        self
          .raw()
          .supersequences(dim_super, parent.raw())
          .map(move |sup| self.complex.skeleton(dim_super).get_by_simplex(&sup))
      })
      .collect()
  }
```

Furthermore these functions always directly access the `IndexSet` and retrieve the corresponding index of the simplex and construct a new `SimplexHandle` out of it, such that we can easily apply subsequent method calls on the returned objects.

## 2.4.2 Boundary Operator

We previously defined the boundary operator $\partial$ for individual simplices. This concept extends naturally to the entire complex. For each dimension $k$, the boundary operator $\partial_k$ maps the $k$-skeleton to the $(k-1)$-skeleton. It can be represented as a linear operator, specifically a matrix $\mathbf{D}_k$, often called the **incidence matrix** [6].

The matrix $\mathbf{D}_k$ has dimensions $N_{k-1} \times N_k$, where $N_j$ is the number of $j$-simplices in the complex. The entry $(\mathbf{D}_k)_{ij}$ represents the signed incidence relation between the $i$-th $(k-1)$-simplex and the $j$-th $k$-simplex. It is $+1$ or $-1$ if the $i$-th $(k-1)$-simplex is a facet of the $j$-th $k$-simplex (with the sign determined by the relative orientation from the boundary definition), and $0$ otherwise.

$$\mathbf{D}_k \in \{-1, 0, +1\}^{N_{k-1} \times N_k} \tag{19}$$

The following function computes this sparse incidence matrix for a given dimension `dim` ($k$).

```rust
impl Complex {
  pub fn boundary_operator(&self, dim: Dim) -> CooMatrix {
    let sups = &self.skeleton(dim);

    if dim == 0 {
      return CooMatrix::zeros(0, sups.len());
    }

    let subs = &self.skeleton(dim - 1);
    let mut mat = CooMatrix::zeros(subs.len(), sups.len());
    for (isup, sup) in sups.handle_iter().enumerate() {
      let sup_boundary = sup.boundary();
      for sub in sup_boundary {
        let sign = sub.sign.as_f64();
        let isub = subs.handle_by_simplex(&sub.simplex).kidx();
        mat.push(isub, isup, sign);
      }
    }
    mat
  }
}
```

## 2.5 Simplicial Geometry

We have now successfully developed the topological structure of our mesh, by combining many abstract simplices into skeletons and collecting all of these skeletons together.

What we are still missing in our mesh data structure now, is any geometry. The geometry is missing, since we only store abstract simplices and not something like coordinate simplices.

This was purposefully done, because we want to separate the topology from the geometry. This allows us to switch between a coordinate-based embedded geometry and a coordinate-free intrinsic geometry based on a Riemannian metric.

## 2.6 Coordinate-Based Simplicial Euclidean Geometry

Let us first quickly look at the familiar coordinate-based euclidean geometry, that relies on an embedding. It's an extrinsic description of the geometry of the manifold from the perspective of the ambient space. All we need is to know the coordinates of all vertices in the mesh.

```rust
#[derive(Debug, Clone)]
pub struct MeshVertexCoords {
  coord_matrix: na::DMatrix<f64>,
}
impl MeshVertexCoords {
  pub fn dim(&self) -> Dim { self.coord_matrix.nrows() }
  pub fn nvertices(&self) -> usize { self.coord_matrix.ncols() }
  pub fn coord(&self, ivertex: VertexIdx) -> CoordRef { self.coord_matrix.column(ivertex) }
}
```

We once again store the coordinates of the vertices in the column of a matrix, just as we did for the `SimplexCoords` struct, but here we store the coordinates of all the vertices in the mesh, so this is usually a really wide matrix. The `dim` function here references the dimension of the ambient space. It is different from the topology dimension in general.

Here we witness another benefit of separating topology and geometry, which should be done even when there are not multiple geometry representations supported: We avoid any redundancy in storing the vertex coordinates. For every vertex we store its coordinate exactly once. This is contrast to using a list of `SimplexCoords`, for which there would have been many duplicate coordinates, since the vertices are shared by many simplices. So separating topology and geometry is always very natural even in the case of the typical coordinate-based geometry.

### 2.6.1 Coordinate Function Functors & Barycentric Quadrature

Before differential geometry, calculus was typically performed on Euclidean space $\mathbb{R}^n$ using global coordinates, rather than on abstract manifolds [5]. A point $p \in \mathbb{R}^n$ is identified with its coordinate vector $x = p$. Functions $f : \Omega \to \mathbb{R}$ are often defined via evaluation rules based on these coordinates, like $f(x) = \sin(x^1)$. This representation, relying on point evaluation, is common in numerical methods. In programming, objects providing such an evaluation capability are often called **functors**.

On general manifolds, global coordinates are usually unavailable [5]. However, if an embedding into an ambient space (RR^N) exists, we can work with functions defined using these ambient coordinates.

A primary application for point-evaluable functions in FEM is numerical integration using **quadrature rules** [1]. Since analytic integration over complex domains or of complicated functions is often infeasible, we approximate integrals numerically.

Quadrature rules for simplices are typically defined on the reference simplex $\hat{\sigma}^n$. A rule consists of a set of quadrature points $q_i$, which are given in local coordinates) and corresponding normalized weights $w_i$, that sum to one. Our implementation uses a `SimplexQuadRule` struct which stores these values.

```
/// A quadrature rule defined on the reference simplex.
pub struct SimplexQuadRule {
  /// Points in local coordinates (reduced barycentric).
  points: na::DMatrix<f64>,
  /// Normalized weights that sum to 1.
  weights: na::DVector<f64>,
}
```

The integral of a function $f$ defined in local coordinates over any simplex $\sigma$ is approximated as:

$$\int_{\hat{\sigma}^n} f \mathrm{vol}_g \approx |\sigma^n| \sum_i w_i f(q_i) \tag{20}$$

Note that here we have factored the volume of the simplex out of the weights, which allows us to store normalized weights that sum to one. Then we can integrate over any simplex, by providing it's volume.

This volume factor $|\sigma|$ accounts for the volume change via the determinant of the Jacobian of the affine map $\varphi_\sigma : \sigma^n \to \sigma$, since $|\sigma| = |\det(\mathrm{D}\varphi_\sigma)||\hat{\sigma}^n|$

```
impl SimplexQuadRule {
  /// Uses a local coordinate function `f`.
  pub fn integrate_local<F>(&self, f: &F, vol: f64) -> f64
  where
    F: Fn(CoordRef) -> f64,
  {
    let mut integral = 0.0;
    for i in 0..self.npoints() {
      integral += self.weights[i] * f(self.points.column(i));
    }
    vol * integral
  }
}
```

To integrate a function defined in ambient coordinates over the simplex, we need to pullback the function.

```
/// Uses a global coordinate function `f`.
pub fn integrate_coord<F>(&self, f: &F, coords: &SimplexCoords) -> f64
where
  F: Fn(CoordRef) -> f64,
{
```

```
    self.integrate_local(
      &|local_coord| f(coords.local2global(local_coord).as_view()),
      coords.vol(),
    )
}
```

To integrate a global function $f$ defined in ambient coordinates over the entire mesh $\mathcal{M} = \cup \, \sigma_j$, we sum the integrals over all cells:

```
/// Uses a global coordinate function `f`.
pub fn integrate_mesh<F>(&self, f: &F, complex: &Complex, coords: &MeshCoords) -> f64
where
  F: Fn(CoordRef) -> f64,
{
  let mut integral = 0.0;
  for cell in complex.cells().iter() {
    let cell_coords = SimplexCoords::from_simplex_and_coords(cell, coords);
    integral += self.integrate_coord(f, &cell_coords);
  }
  integral
}
```

The library provides several standard quadrature rules.

The most important of which is the barycentric quadrature rule. It generalizes trivially to arbitrary dimensional $n$-simplices. This rule has polynomial exactness degree 1, meaning it integrates affine linear functions exactly. This is sufficient for most applications in 1st order FEEC, since it constitutes an admissible variational crime. [1]

```
/// Integrates 1st order affine linear functions exactly.
pub fn barycentric(dim: Dim) -> Self {
  let barycenter = barycenter_local(dim);
  let points = Matrix::from_columns(&[barycenter]);
  let weight = 1.0;
  let weights = Vector::from_element(1, weight);
  Self { points, weights }
}
```

We also have hard-coded order 3 rules for 1D, 2D and 3D. They are accurate for polynomials of degree 2. This will be crucial for measuring the convergence of our FE solutions to the exact solution, for which we need to use a quadrature rule, where the quadrature error doesn't dominate.

```
pub fn order3(dim: Dim) -> Self {
  match dim {
    0 => Self::dim0(),
    1 => Self::dim1_order3(),
    2 => Self::dim2_order3(),
    3 => Self::dim3_order3(),
    _ => unimplemented!("No order 3 Quadrature available for dim {dim}."),
  }
}

/// Simpsons Rule
pub fn dim1_order3() -> Self {
  let points = na::dmatrix![0.0, 0.5, 1.0];
  let weights = na::dvector![1.0 / 6.0, 4.0 / 6.0, 1.0 / 6.0];
  Self { points, weights }
}
pub fn dim2_order3() -> Self {
  let points = na::dmatrix![
    1.0/3.0, 1.0/5.0, 3.0/5.0, 1.0/5.0;
    1.0/3.0, 1.0/5.0, 1.0/5.0, 3.0/5.0;
  ];
  let weights = na::dvector![-27.0 / 48.0, 25.0 / 48.0, 25.0 / 48.0, 25.0 / 48.0];
  Self { points, weights }
}
pub fn dim3_order3() -> Self {
```

```
    // omitted...
}
```

## 2.7 Metric-Based Riemannian Geometry

The coordinate-based Euclidean geometry we've seen so far, is what is commonly used in almost all FEM implementations. In our implementation we go one step further and abstract away the coordinates of the manifold and instead make use of coordinate-free Riemannian geometry [5] for all of our FE algorithms. All FE algorithms only depend on this geometry representation and cannot operate directly on the coordinate-based geometry. Instead one should always derive a coordinate-free representation from the coordinate-based one. Most of the time one starts with a coordinate-based representation that has been constructed by some mesh generator like gmsh [21] and then one computes the intrinsic geometry and forgets about the coordinates. Our library supports this exactly this functionality.

### 2.7.1 Riemannian Metric

Riemannian geometry is an intrinsic description of the manifold, that doesn't need an ambient space at all. It relies purely on a structure over the manifold called a **Riemannian metric** $g$ [5].

It is a continuous function over the whole manifold, which at each point $p$ gives us an inner product $g_p : T_pM \times T_pM \to \mathbb{R}^+$ on the tangent space $T_pM$ at this point $p$. It is the analog to the standard euclidean inner product (dot product) in euclidean geometry. The inner product on tangent vectors allows one to measure lengths $\|v\|_g = \sqrt{g(v,v)}$ angles $\varphi(v,w) = \arccos\left(\frac{g_p(v,w)}{\|v\|_g\|w\|_g}\right)$. While euclidean space is flat and the inner product is the same everywhere, a manifold is curved in general and therefore the inner product changes from point to point, reflecting the changing geometry.

Given a basis $\frac{\partial}{\partial x^1}, ..., \frac{\partial}{\partial x^n}$ of the tangent space $T_pM$ at a point $p$, induced by a chart map $\varphi : p \in U \subseteq M \mapsto (x_1, ..., x_n)$, the inner product $g_p$ can be represented in this basis using components $\left(g_p\right)_{ij} \in \mathbb{R}$. This is done by plugging in all combinations of basis vectors into the two arguments of the bilinear form.

$$\left(g_p\right)_{ij} = g_p\left(\left.\frac{\partial}{\partial x^i}\right|_p, \left.\frac{\partial}{\partial x^j}\right|_p\right) \tag{21}$$

We can collect all these components into a matrix $\mathbf{G} \in \mathbb{R}^{n \times n}$

$$\mathbf{G} = \left[g_{ij}\right]_{i,j=1}^{n \times n} \tag{22}$$

This is called a Gram matrix or Gramian and is the discretization of a inner product of a linear space, given a basis. This matrix doesn't represent a linear map, which would be a $(1,1)$-tensor, but instead a bilinear form, which is a $(0,2)$-tensor. In the context of Riemannian geometry this is called a **metric tensor** [5].

The inverse metric $g_p^{-1}$ at a point $p$ provides an inner product $g_p^{-1} : T_p^*M \times T_p^*M \to \mathbb{R}^+$ on the cotangent space $T_p^*M$. It can be obtained by computing the inverse Gramian matrix $\mathbf{G}^{-1}$, which is then a new Gramian matrix representing the inner product on the dual basis of covectors.

$$\mathbf{G}^{-1} = \left[g(\mathrm{d}x^i, \mathrm{d}x^j)\right]_{i,j=1}^{n \times n} \tag{23}$$

The inverse metric is very important for us, since differential forms are covariant tensors, therefore they are measured by the inverse metric tensor [5].

Computing the inverse is numerically unstable and instead it would be better to rely on matrix factorization to do computation involving the inverse metric [23]. However this quickly becomes intractable. For this reason we chose here to rely on the directly computed inverse matrix nonetheless.

We introduce a struct to represent the Riemannian metric at a particular point as the Gramian matrix.

```
/// A Gram Matrix represents an inner product expressed in a basis.
#[derive(Debug, Clone)]
pub struct Gramian {
  /// S.P.D. matrix
  matrix: na::DMatrix<f64>,
}
```

```rust
impl Gramian {
  pub fn try_new(matrix: na::DMatrix<f64>) -> Option<Self> {
    matrix.is_spd().then_some(Self { matrix })
  }
  pub fn new(matrix: na::DMatrix<f64>) -> Self {
    Self::try_new(matrix).expect("Matrix must be s.p.d.")
  }
  pub fn new_unchecked(matrix: na::DMatrix<f64>) -> Self {
    if cfg!(debug_assertions) {
      Self::new(matrix)
    } else {
      Self { matrix }
    }
  }
  pub fn from_euclidean_vectors(vectors: na::DMatrix<f64>) -> Self {
    assert!(vectors.is_full_rank(1e-9), "Matrix must be full rank.");
    let matrix = vectors.transpose() * vectors;
    Self::new_unchecked(matrix)
  }
  /// Orthonormal euclidean metric.
  pub fn standard(dim: Dim) -> Self {
    let matrix = na::DMatrix::identity(dim, dim);
    Self::new_unchecked(matrix)
  }

  pub fn matrix(&self) -> &na::DMatrix<f64> { &self.matrix }
  pub fn dim(&self) -> Dim { self.matrix.nrows() }
  pub fn det(&self) -> f64 { self.matrix.determinant() }
  pub fn det_sqrt(&self) -> f64 { self.det().sqrt() }

  pub fn inverse(self) -> Self {
    let matrix = self
      .matrix
      .try_inverse()
      .expect("Symmetric Positive Definite is always invertible.");
    Self::new_unchecked(matrix)
  }
}

/// Inner product functionality directly on the basis.
impl Gramian {
  pub fn basis_inner(&self, i: usize, j: usize) -> f64 { self.matrix[(i, j)] }
  pub fn basis_norm_sq(&self, i: usize) -> f64 { self.basis_inner(i, i) }
  pub fn basis_norm(&self, i: usize) -> f64 { self.basis_norm_sq(i).sqrt() }
  pub fn basis_angle_cos(&self, i: usize, j: usize) -> f64 {
    self.basis_inner(i, j) / self.basis_norm(i) / self.basis_norm(j)
  }
  pub fn basis_angle(&self, i: usize, j: usize) -> f64 { self.basis_angle_cos(i, j).acos() }
}
impl std::ops::Index<(usize, usize)> for Gramian {
  type Output = f64;
  fn index(&self, (i, j): (usize, usize)) -> &Self::Output {
    &self.matrix[(i, j)]
  }
}

/// Inner product functionality directly on any element.
impl Gramian {
  pub fn inner(&self, v: &na::DVector<f64>, w: &na::DVector<f64>) -> f64 {
    (v.transpose() * self.matrix() * w).x
  }
  pub fn inner_mat(&self, v: &na::DMatrix<f64>, w: &na::DMatrix<f64>) -> na::DMatrix<f64> {
    v.transpose() * self.matrix() * w
  }
  pub fn norm_sq(&self, v: &na::DVector<f64>) -> f64 {
```

```rust
    self.inner(v, v)
  }
  pub fn norm_sq_mat(&self, v: &na::DMatrix<f64>) -> na::DMatrix<f64> {
    self.inner_mat(v, v)
  }
  pub fn norm(&self, v: &na::DVector<f64>) -> f64 {
    self.inner(v, v).sqrt()
  }
  pub fn norm_mat(&self, v: &na::DMatrix<f64>) -> na::DMatrix<f64> {
    self.inner_mat(v, v).map(|v| v.sqrt())
  }
  pub fn angle_cos(&self, v: &na::DVector<f64>, w: &na::DVector<f64>) -> f64 {
    self.inner(v, w) / self.norm(v) / self.norm(w)
  }
  pub fn angle(&self, v: &na::DVector<f64>, w: &na::DVector<f64>) -> f64 {
    self.angle_cos(v, w).acos()
  }
}
```

The dot product is the standard inner product on flat Euclidean space. The standard basis vectors are orthonormal w.r.t. this inner product. Therefore the gram matrix (and it's inverse) are just identity matrices.

```rust
/// Orthonormal flat euclidean metric.
pub fn standard(dim: Dim) -> Self {
  let identity = na::DMatrix::identity(dim, dim);
  let metric_tensor = identity.clone();
  let inverse_metric_tensor = identity;
  Self { metric_tensor, inverse_metric_tensor, }
}
```

### 2.7.2 Deriving the Metric from an Immersion

One can easily derive the Riemannian metric from an immersion $f : M \to \mathbb{R}^N$ into an ambient space $\mathbb{R}^N$ [5]. It's differential is a function $\mathrm{d}f_p : T_p M \to T_{\mathbf{p}} \mathbb{R}^N$, also called the push-forward and tells us how our intrinsic tangential vectors are being placed into the ambient space, giving them an extrinsic geometry.

This immersion then induces a metric, that describes the same geometry intrinsically. For this we just take the standard euclidean inner product of our immersed tangent vectors. This then inherits the extrinsic ambient geometry and represents it intrinsically.

$$g_p(u, v) = \mathrm{d}f_p(u) \cdot \mathrm{d}f_p(v) \tag{24}$$

Computationally this differential $\mathrm{d}f_p$ can be represented, given a basis, since it is a linear map, by a Jacobi Matrix $\mathbf{J}$. The metric is the then the Gramian matrix of the Jacobian.

$$\mathbf{G} = \mathbf{J}^\top \mathbf{J} \tag{25}$$

The Jacobian has as columns our immersed basis tangent vectors, therefore really we just need these to compute a metric.

```rust
pub fn from_euclidean_vectors(vectors: na::DMatrix<f64>) -> Self
```

## 2.8 Simplicial Riemannian Geometry & Regge Calculus

We have now discussed the of Riemannian geometry in general. Now we want to focus on the special case of simplicial geometry, where the Riemannian metric is defined on our mesh.

We have seen with our coordinate simplices that our geometry is piecewise-flat over the cells. This means that our metric is constant over each cell and changes only from one cell to another.

This piecewise-constant metric over the simplicial mesh can be derived either from vertex coordinates (relying on an embedding) or intrinsically from the edge lengths of the simplex by a method inspired by **Regge calculus** [9], a theory for numerical general relativity that is about producing simplicial approximations of spacetimes that are solutions to the Einstein field equation.

### 2.8.1 Deriving the Metric from Coordinate Simplices

Our coordinate simplices are an immersion of an abstract simplex and as such, we can compute the corresponding constant metric tensor on it. We've seen that the local metric tensor is exactly the Gramian of the Jacobian of the immersion. The Jacobian of our affine-linear transformation, is the matrix $\mathbf{E}$, that has the spanning vectors as columns. From a different perspective the spanning vectors just constitute our chosen basis of the tangent space and therefore it's Gramian w.r.t. the Euclidean inner product is the metric tensor.

$$\mathbf{G} = \mathbf{E}^\top \mathbf{E}$$
$$\mathbf{G}_{ij} = e_i \cdot e_j \tag{26}$$

```rust
impl SimplexCoords {
  pub fn metric_tensor(&self) -> Gramian {
    Gramian::from_euclidean_vectors(self.spanning_vectors())
  }
}
```

### 2.8.2 Simplex Edge Lengths

We have seen how to derive the metric tensor from an embedding, but of course the metric is independent of the specific coordinates. The metric is invariant under isometric transformations, such as translations, rotations and reflections [5]. This begs the question, what the minimal geometric information necessary is to derive the metric. It turns out that, while vertex coordinates are over-specified, edge lengths in contrast are exactly the required information. Edge lengths are invariant under isometric transformations. This is known from **Regge Calculus** [9].

Instead of giving all vertices a global coordinate, as one would do in extrinsic geometry, we just give each edge in the mesh a positive length. Just knowing the lengths doesn't tell you the positioning of the mesh in an ambient space but it's enough to give the whole mesh it's piecewise-flat geometry.

Mathematically this is just a function on the edges to the positive real numbers.

$$d : \Delta_1(\mathcal{M}) \to \mathbb{R}^+$$
$$d_{ij} = d\big([v_i, v_j]\big) \tag{27}$$

that gives each edge $e \in \Delta_1(\mathcal{M})$ a positive length $l_e \in \mathbb{R}^+$.
We denote the edge length between vertices $i$ and $j$ by $d_{ij}$.

```rust
#[derive(Debug, Clone)]
pub struct SimplexLengths {
  lengths: na::DVector<f64>,
  dim: Dim,
}

pub fn new(lengths: na::DVector<f64>, dim: Dim) -> Self {
  assert_eq!(lengths.len(), nedges(dim), "Wrong number of edges.");
  let this = Self { lengths, dim };
  assert!(
    this.is_coordinate_realizable(),
    "Simplex must be coordinate realizable."
  );
  this
}
pub fn new_unchecked(lengths: na::DVector<f64>, dim: Dim) -> Self {
  if cfg!(debug_assertions) {
    Self::new(lengths, dim)
  } else {
    Self { lengths, dim }
  }
}
pub fn standard(dim: Dim) -> SimplexLengths {
  let nedges = nedges(dim);
  let lengths: Vec<f64> = (0..dim)
    .map(|_| 1.0)
    .chain((dim..nedges).map(|_| SQRT_2))
```

```
    .collect();

  Self::new_unchecked(lengths.into(), dim)
}
pub fn from_coords(coords: &SimplexCoords) -> Self {
  let dim = coords.dim_intrinsic();
  let lengths = coords.edges().map(|e| e.vol()).collect_vec().into();
  // SAFETY: Edge lengths stem from a realization already.
  Self::new_unchecked(lengths, dim)
}
```

The edge lengths of a simplex and the metric tensor on it both fully define the geometry uniquely [9]. These two geometry representations are completely equivalent. This means one can be derived from the other.

We can derive the edge lengths, from the metric tensor Gramian.

$$d_{ij} = \sqrt{\mathbf{G}_{ii} + \mathbf{G}_{jj} - 2\mathbf{G}_{ij}} \qquad (28)$$

```
pub fn from_metric_tensor(metric: &Gramian) -> Self {
  let dim = metric.dim();
  let length = |i, j| {
    (metric.basis_inner(i, i) + metric.basis_inner(j, j) - 2.0 * metric.basis_inner(i, j)).sqrt()
  };

  let mut lengths = na::DVector::zeros(nedges(dim));
  let mut iedge = 0;
  for i in 0..dim {
    for j in i..dim {
      lengths[iedge] = length(i, j);
      iedge += 1;
    }
  }

  Self::new(lengths, dim)
}
```

We can derive the metric, from the edge lengths by using the law of cosines.

$$\mathbf{G}_{ij} = \frac{1}{2}\left(d_{0i}^2 + d_{0j}^2 - d_{ij}^2\right) \qquad (29)$$

```
pub fn to_metric_tensor(&self) -> Gramian {
  let mut metric = Matrix::zeros(self.dim(), self.dim());
  for i in 0..self.dim() {
    metric[(i, i)] = self[i].powi(2);
  }
  for i in 0..self.dim() {
    for j in (i + 1)..self.dim() {
      let l0i = self[i];
      let l0j = self[j];

      let vi = i + 1;
      let vj = j + 1;
      let eij = lex_rank(&[vi, vj], self.nvertices());
      let lij = self[eij];

      let val = 0.5 * (l0i.powi(2) + l0j.powi(2) - lij.powi(2));

      metric[(i, j)] = val;
      metric[(j, i)] = val;
    }
  }
  Gramian::new(metric)
}
```

### 2.8.3 Realizability Conditions

While edge lengths provide a coordinate-free description of a simplex's intrinsic geometry, not just any assignment of positive numbers to the edges of an abstract simplex constitutes a valid Euclidean geometry. The assigned lengths must satisfy certain consistency requirements, known collectively as **realizability conditions** [24]. These conditions ensure that the abstract simplex, endowed with these edge lengths, could actually be embedded isometrically as a flat simplex in some Euclidean space $\mathbb{R}^N$. In essence, only edge lengths for which a corresponding `SimplexCoords` object could exist could be considered valid for defining a Euclidean simplex geometry.

A fundamental necessary condition stems from the geometry of triangles. For every 2-dimensional face (triangle) $\sigma = [i, j, k]$ within the mesh, the assigned edge lengths $d_{ij}$, $d_{jk}$, and $d_{ik}$ must satisfy the standard **triangle inequalities**:

$$d_{ij} + d_{jk} \geq d_{ik}$$
$$d_{jk} + d_{ik} \geq d_{ij} \tag{30}$$
$$d_{ik} + d_{ij} \geq d_{jk}$$

If these inequalities (particularly the strict versions) are violated for any triangle, the metric tensor $\mathbf{G}$ derived from these lengths via the law of cosines will not be positive-definite. This would imply a degenerate or pseudo-Riemannian metric rather than the proper Riemannian metric associated with Euclidean geometry. If the edge lengths are derived from an actual coordinate embedding (`SimplexCoords::from_coords`), the triangle inequalities are automatically satisfied.

While necessary, the triangle inequalities alone are not sufficient for dimensions $n > 2$. A more comprehensive check involves the squared distances between all pairs of vertices. We can assemble these into the **Euclidean distance matrix** (EDM) $\mathbf{A}$ [24], a symmetric matrix with zeros on the diagonal:

$$\mathbf{A} = \begin{bmatrix} 0 & d_{12}^2 & d_{13}^2 & \cdots & d_{1n}^2 \\ d_{21}^2 & 0 & d_{23}^2 & \cdots & d_{2n}^2 \\ d_{31}^2 & d_{32}^2 & 0 & \cdots & d_{3n}^2 \\ \vdots & \vdots & \vdots & \ddots & \vdots \\ d_{n1}^2 & d_{n2}^2 & d_{n3}^2 & \cdots & 0 \end{bmatrix} \tag{31}$$

where $d_{ij}$ is the length of the edge between vertex $i$ and vertex $j$.

```
pub fn distance_matrix(&self) -> na::DMatrix<f64> {
  let mut mat = na::DMatrix::zeros(self.nvertices(), self.nvertices());

  let mut idx = 0;
  for i in 0..self.nvertices() {
    for j in (i + 1)..self.nvertices() {
      let dist_sqr = self.lengths[idx].powi(2);
      mat[(i, j)] = dist_sqr;
      mat[(j, i)] = dist_sqr;
      idx += 1;
    }
  }
  mat
}
```

Building upon the EDM, the **Cayley-Menger matrix** [24] is constructed by bordering the EDM with a row and column of ones:

$$\mathbf{CM} = \begin{bmatrix} 0 & d_{12}^2 & d_{13}^2 & \cdots & d_{1n}^2 & 1 \\ d_{21}^2 & 0 & d_{23}^2 & \cdots & d_{2n}^2 & 1 \\ d_{31}^2 & d_{32}^2 & 0 & \cdots & d_{3n}^2 & 1 \\ \vdots & \vdots & \vdots & \ddots & \vdots & \vdots \\ d_{n1}^2 & d_{n2}^2 & d_{n3}^2 & \cdots & 0 & 1 \\ 1 & 1 & 1 & \cdots & 1 & 0 \end{bmatrix} \tag{32}$$

```
pub fn cayley_menger_matrix(&self) -> na::DMatrix<f64> {
  let mut mat = self.distance_matrix();
  mat = mat.insert_row(self.nvertices(), 1.0);
  mat = mat.insert_column(self.nvertices(), 1.0);
```

```rust
    mat[(self.nvertices(), self.nvertices())] = 0.0;
    mat
}
```

The determinant of this matrix, scaled by a dimension-dependent factor, is the **Cayley-Menger determinant** $cm$ [24]:

$$cm = \frac{(-1)^{n+1}}{(n!)^2 2^n} \det(\mathbf{CM}) \tag{33}$$

```rust
impl SimplexLengths {
  pub fn cayley_menger_det(&self) -> f64 {
    cayley_menger_factor(self.dim()) * self.cayley_menger_matrix().determinant()
  }
}
pub fn cayley_menger_factor(dim: Dim) -> f64 {
  (-1.0f64).powi(dim as i32 + 1) / factorial(dim).pow(2) as f64 / 2f64.powi(dim as i32)
}
```

A fundamental result states that a set of edge lengths for an $n$-simplex is realizable in Euclidean space $\mathbb{R}^N$ (for $N \geq n$) if and only if the Cayley-Menger determinant is non-negative: $cm \geq 0$ [24]. A strictly positive determinant indicates realizability in exactly $n$ dimensions (non-degenerate), while a zero determinant implies the simplex is degenerate and lies within an lower dimensional affine subspace.

```rust
pub fn is_coordinate_realizable(&self) -> bool {
  self.cayley_menger_det() >= 0
}
```

Furthermore, if the simplex is realizable ($cm \geq 0$), its $n$-dimensional volume is directly related to the Cayley-Menger determinant:

$$\mathrm{vol} = \sqrt{cm} \tag{34}$$

This provides a way to compute the volume directly from edge lengths, without needing explicit coordinates or the metric tensor.

```rust
pub fn vol(&self) -> f64 {
  self.cayley_menger_det().sqrt()
}
```

### 2.8.4 Global Geometry

Computationally we represent the edge lengths in a single struct that has all lengths stored continuously in memory in a nalgebra vector [14].

```rust
pub struct MeshEdgeLengths {
  vector: na::DVector<f64>,
}
```

Our topological simplicial complex struct gives a global numbering to our edges, which then gives us the indices into this nalgebra vector.

Our topological simplicial manifold together with these edge lengths gives us a simplicial Riemannian manifold.

Our FE algorithms than usually take two arguments.

```rust
fn fe_algorithm(topology: &Complex, geometry: &MeshEdgeLengths)
```

We can derive a MeshEdgeLengths struct from a MeshVertexCoords struct.

```rust
impl MeshVertexCoords {
  pub fn to_edge_lengths(&self, topology: &Complex) -> MeshEdgeLengths {
    let edges = topology.edges();
    let mut edge_lengths = na::DVector::zeros(edges.len());
    for (iedge, edge) in edges.set_iter().enumerate() {
      let [vi, vj] = edge.clone().try_into().unwrap();
      let length = (self.coord(vj) - self.coord(vi)).norm();
      edge_lengths[iedge] = length;
    }
```

```
    MeshEdgeLengths::new(edge_lengths)
  }
}
```

## 2.9 Higher Order Geometry

The original PDE problem, before discretization, is posed on a smooth manifold, which we then discretized in the form of a mesh. This smooth manifold has a non-zero curvature everywhere in general [5]. This is in contrast to the simplex geometry we've chosen here, that approximates the geometry, by a **piecewise-flat geometry** over the cells. Each cell is a flat simplex that has no curvature change inside of it.

This manifests as the fact, that for a coordinate-based representation of the geometry, the cell is just the convex hull of the vertex coordinates. And for the metric-based representation, we have a **piecewise-constant metric** over each cell.

This is a piecewise-linear 1st order approximation of the geometry. But this is not only possible approximation. Higher-order mesh elements, such as quadratic or cubic elements allow for higher accuracy approximations of the mesh [1]. In general any order polynomial elements can be used. We will restrict ourselves completely to first-order elements in this thesis. This approximation is sufficient for us, since it represents an **admissible geometric variational crime** [11]: The order of our FE method coincides with the order of our mesh geometry; both are linear 1st order. This approximation doesn't affect the order of convergence in a negative way, and therefore is admissible [11].

## 2.10 Mesh Generation and Loading

### 2.10.1 Tensor-Product Domain Meshing

Formoniq features a meshing algorithm for arbitrary dimensional tensor-product domains. These domains are $n$-dimensional Cartesian products $[0,1]^n$ of the unit interval $[0,1]$. The simplicial skeleton will be computed based on a Cartesian grid that subdivides the domain into $l^n$ many $n$-cubes, which are generalizations of squares and cubes. Here $l$ is the number of subdivisions per axis. To obtain a simplicial skeleton, we need to split each $n$-cube into non-overlapping $n$-simplices that make up its volume. In 2D it's very natural to split a square into two triangles of equal volume. This can be generalized to higher dimensions. The trivial triangulation of a $n$-cube into $n!$ simplices is based on the $n!$ many permutations of the $n$ coordinate axes. [25]

The $n$-cube has $2^n$ vertices, which can all be identified using multiindices

$$V = \{0,1\}^n = \left\{ (i_1, ..., i_n) \,\middle|\, i_j \in \{0,1\} \right\} \tag{35}$$

All $n!$ simplices will be based on this vertex base set. To generate the list of vertices of the simplex, we start at the origin vertex $v_0 = 0 = (0)^n$. From there we walk along axis directions from vertex to vertex. For this we consider all $n!$ permutations of the basis directions $e_1, ..., e_n$. A permutation $\sigma$ tells in which axis direction we need to walk next. This gives us the vertices $v_0, ..., v_n$ that forms a simplex.

$$v_k = v_0 + \sum_{i=1}^{k} \boldsymbol{e}_{\sigma(i)} \tag{36}$$

The algorithm in Rust looks like:

```rust
pub fn compute_cell_skeleton(&self) -> Skeleton {
  let nboxes = self.ncells();
  let nboxes_axis = self.ncells_axis();

  let dim = self.dim();
  let nsimplices = factorial(dim) * nboxes;
  let mut simplices: Vec<SortedSimplex> = Vec::with_capacity(nsimplices);

  // iterate through all boxes that make up the mesh
  for ibox in 0..nboxes {
    let cube_icart = linear_index2cartesian_index(ibox, nboxes_axis, self.dim());

    let vertex_icart_origin = cube_icart;
    let ivertex_origin =
```

```
        cartesian_index2linear_index(vertex_icart_origin.clone(), self.nvertices_axis());

    let basisdirs = IndexSet::increasing(dim);

    // Construct all $d!$ simplexes that make up the current box.
    // Each permutation of the basis directions (dimensions) gives rise to one simplex.
    let cube_simplices = basisdirs.permutations().map(|basisdirs| {
      // Construct simplex by adding all shifted vertices.
      let mut simplex = vec![ivertex_origin];

      // Add every shift (according to permutation) to vertex iteratively.
      // Every shift step gives us one vertex.
      let mut vertex_icart = vertex_icart_origin.clone();
      for basisdir in basisdirs.set.iter() {
        vertex_icart[basisdir] += 1;

        let ivertex = cartesian_index2linear_index(vertex_icart.clone(), self.nvertices_axis());
        simplex.push(ivertex);
      }

      Simplex::from(simplex).assume_sorted()
    });

    simplices.extend(cube_simplices);
  }

  Skeleton::new(simplices)
}
```

We can note here, that the computational complexity of this algorithm, grows extremely fast in the dimension $n$. We have a factorial scaling $\mathcal{O}(n!)$ (worse than exponential scaling $\mathcal{O}(e^n)$) for splitting the cube into simplices. Given $l$ subdivisions per dimensions, we have $l^n$ cubes. So the overall computational complexity is dominated by $\mathcal{O}(l^n n!)$, a terrible result, due to the curse of dimensionality. The memory usage is dictated by the same scaling law.

### 2.10.2 Mesh Import

The formoniq manifold crate supports multiple mesh formats for importing.

Formoniq can read gmsh `.msh` files [21] and turn them into a simplicial complex that we can work on. This is thanks to the `gmshio` crate [22].

```
pub fn gmsh2coord_complex(bytes: &[u8]) -> (Complex, MeshCoords)
```

Formoniq also supports reading and writing `.obj` files, which are very common in computer graphics software, such as blender. We implemented some very basic custom readers and writers, since the format is very simple.

```
pub fn obj2coord_complex(obj_string: &str) -> (Complex, MeshCoords)
```

# Chapter 3

# Exterior Algebra

Exterior algebra is to exterior calculus, what vector algebra is to vector calculus.

In vector calculus we have vector fields $v$, which are functions $v : p \in \Omega \mapsto v_p$ over the manifold $\Omega$ that at each point $p$ have a constant vector $v_p \in T_p M$ as value.

In exterior calculus we have differential forms $\omega$, which are functions $\omega : p \in \Omega \mapsto \omega_p$ over the manifold $\Omega$ that at each point $p$ have a constant **multiform** $\omega_p \in \bigwedge(T_p^* M)$ as value [5], [8].

If one were to implement something related to vector calculus it is of course crucial to be able to represent vectors in the program. This is usually the job of a basic linear algebra library such as Eigen in C++ and nalgebra in Rust [14]. Since we want to implement FEEC, which uses exterior calculus, it is crucial, that we are able to represent multiforms in our program. For this there aren't any established libraries. So we do this ourselves and develop a small module.

## 3.1  Exterior Algebra of Multiforms

In general an exterior algebra $\bigwedge(V)$ is a construction over any linear space $V$ [5]. In this section we want to quickly look at the specific linear space we are dealing with when modelling multiforms as element of an exterior algebra. But our implementation would work for any finite-dimensional real linear space $V$ with a given basis.

In our particular case we have the exterior algebra of alternating multilinear forms $\bigwedge(T_p^* M)$ [5]. Here the linear space $V$ is the cotangent space $T_p^* M$ of the manifold $M$ at a point $p \in M$. It's the dual space $(T_p M)^*$ of the tangent space $T_p M$. The elements of the cotangent space are covectors $a \in T_p^* M$, which are linear functionals $a : T_p M \to \mathbb{R}$ on the tangent space. The tangent space $T_p M$ has the standard basis $\left\{ \frac{\partial}{\partial x^i} \right\}_{i=1}^n$ induced by some chart map $\varphi : p \in U \subseteq M \mapsto (x_1, ..., x_n)$. This gives rise to a dual basis $\left\{ \mathrm{d}x^i \right\}_{i=1}^n$ of the cotangent space, defined by $\mathrm{d}x^i \left( \frac{\partial}{\partial x^j} \right) = \delta_j^i$ [5].

There is a related space, called the space of multivectors $\bigwedge(T_p M)$, which is the exterior algebra over the tangent space, instead of the cotangent space. The space of multivectors and multiforms are dual to each other [5].

$$\bigwedge(T_p^* M) \cong \left( \bigwedge(T_p M) \right)^* \tag{37}$$

The space of multivectors only plays a minor role in exterior calculus.

The elements of an exterior algebra are commonly called multivectors, irregardless what the underlying linear space $V$ is. This is confusing when working with multiforms, which are distinct from multivectors. To avoid confusion, we therefore just call the elements of the exterior algebra exterior elements or multielements, just like we say linear space instead of vector space.

## 3.2  The Numerical Exterior Algebra $\bigwedge(\mathbb{R}^n)$

When working with vectors from a finite-dimensional real linear space $V$, then we can always represent them computationally, by choosing a basis $\{e_i\}_{i=1}^n \subseteq V$. This constructs an isomorphism $V \cong \mathbb{R}^n$, where $n = \dim V$. This allows us to work with elements $v \in \mathbb{R}^n$, which have real values $v_i \in \mathbb{R}$ as components, which are the basis coefficients. These real numbers are what we can work with on computers and allow us to do numerical linear algebra [23]. This means that when working with any finite-dimensional real linear space $V$ on a computer we always just use the linear space $\mathbb{R}^n$.

The same idea can be used to computationally work with exterior algebras [8]. By choosing a basis of $V$, we also get an isomorphism on the exterior algebra $\bigwedge(V) \cong \bigwedge(\mathbb{R}^n)$. Therefore we use $\bigwedge(\mathbb{R}^n)$ in our implementation.

For our space of multiforms, we will be using the standard cotangent basis $\{\mathrm{d}x^i\}_{i=1}^n$.

## 3.3 Representing Exterior Elements

An exterior algebra is a graded algebra [5].

$$\bigwedge(\mathbb{R}^n) = \bigwedge\nolimits^0(\mathbb{R}^n) \bigoplus \cdots \bigoplus \bigwedge\nolimits^n(\mathbb{R}^n) \tag{38}$$

Each element $v \in \bigwedge(\mathbb{R}^n)$ has some particular exterior grade $k \in \{1, ..., n\}$ and therefore lives in a particular exterior power $v \in \bigwedge^k(\mathbb{R}^n)$. We make use of this fact in our implementation, by splitting the representation between these various grades.

```
pub type ExteriorGrade = usize;
```

For representing an element in a particular exterior power $\bigwedge^k(\mathbb{R}^n)$, we use the fact that, it itself is a linear space in it's own right. Due to the combinatorics of the anti-symmetric exterior algebra, we have $\dim \bigwedge^k(\mathbb{R}^n) = \binom{n}{k}$ [5]. This means that by choosing a basis $\{e_I\}$ of this exterior power, we can just use a list of $\binom{n}{k}$ coefficients to represent an exterior element, by using the isomorphism $\bigwedge^k(\mathbb{R}^n) \cong \mathbb{R}^{\binom{n}{k}}$.

```
/// An element of an exterior algebra.
#[derive(Debug, Clone)]
pub struct ExteriorElement {
  coeffs: na::DVector<f64>,
  dim: Dim,
  grade: ExteriorGrade,
}
```

This struct represents an element `self` $\in \bigwedge^k(\mathbb{R}^n)$ with `self.dim` $= n$, `self.grade` $= k$ and `self.coeffs.len()` $= \binom{n}{k}$.

This exterior basis $\{e_I\}_{I \in \mathcal{J}_k^n}$ is different from the basis $\{e_i\}_{i=1}^n$ of the original linear space $V$, but is best subsequently constructed from it. We do this by creating elementary multielements from the exterior product of basis elements [5].

$$e_I = \bigwedge_{j=1}^k e_{I_j} = e_{i_1} \wedge \cdots \wedge e_{i_k} \tag{39}$$

Here $I = (i_1, ..., i_k)$ is a multiindex, in particular a $k$-index, telling us which basis elements to wedge.

Because of the anti-symmetry of the exterior product, there are certain conditions on the multiindices $I$ for $\{e_I\}$ to constitute a meaningful basis. First $I$ must not contain any duplicate indices, because otherwise $e_I = 0$ and second there must not be any permutations of the same index in the basis set, otherwise we have linear dependence of the two elements. We therefore only consider strictly increasing multiindices $I \in \mathcal{J}_k^n$ and denote their set by $\mathcal{J}_k^n = \{(i_1, ..., i_k) \in \mathbb{N}^k \mid 1 \le i_1 < \cdots < i_k \le n\}$. This works for arbitrary dimensions.

The basis also needs to be ordered, such that we can know which coefficient in `self.coeffs` corresponds to which basis. A natural choice here is a lexicographical ordering [8].

Taking in all of this together we for example have as exterior basis for $\bigwedge^2(\mathbb{R}^3)$ the elements $e_1 \wedge e_2, e_1 \wedge e_3, e_2 \wedge e_3$.

## 3.4 Representing Exterior Terms

It is helpful to represent these exterior basis wedges in our program.

```
#[derive(Debug, Clone, PartialEq, Eq, Hash)]
pub struct ExteriorTerm {
  indices: Vec<usize>,
  dim: Dim,
}
impl ExteriorTerm {
  pub fn dim(&self) -> Dim { self.dim }
```

```rust
  pub fn grade(&self) -> ExteriorGrade { self.indices.len() }
}
```

This struct allows for any multiindex, even if they are not strictly increasing. But we are of course able to check whether this is the case or not and then to convert it into a increasing representation plus the permutation sign. We call this representation, canonical.

```rust
pub fn is_basis(&self) -> bool {
  self.is_canonical()
}
pub fn is_canonical(&self) -> bool {
  let Some((sign, canonical)) = self.clone().canonicalized() else {
    return false;
  };
  sign == Sign::Pos && canonical == *self
}
pub fn canonicalized(mut self) -> Option<(Sign, Self)> {
  let sign = sort_signed(&mut self.indices);
  let len = self.indices.len();
  self.indices.dedup();
  if self.indices.len() != len {
    return None;
  }
  Some((sign, self))
}
```

In the case of a strictly increasing term, we can also determine the lexicographical rank of it in the set of all increasing terms. And the other way constructing them from lexicographical rank.

```rust
pub fn lex_rank(&self) -> usize {
  assert!(self.is_canonical(), "Must be canonical.");
  let n = self.dim();
  let k = self.indices.len();

  let mut rank = 0;
  for (i, &index) in self.indices.iter().enumerate() {
    let start = if i == 0 { 0 } else { self.indices[i - 1] + 1 };
    for s in start..index {
      rank += binomial(n - s - 1, k - i - 1);
    }
  }
  rank
}

pub fn from_lex_rank(dim: Dim, grade: ExteriorGrade, mut rank: usize) -> Self {
  let mut indices = Vec::with_capacity(grade);
  let mut start = 0;
  for i in 0..grade {
    let remaining = grade - i;
    for x in start..=(dim - remaining) {
      let c = binomial(dim - x - 1, remaining - 1);
      if rank < c {
        indices.push(x);
        start = x + 1;
        break;
      } else {
        rank -= c;
      }
    }
  }
  Self::new(indices, dim)
}
```

Now that we have this we can implement a useful iterator on our ExteriorElement struct that allows us to iterate through the basis expansion consisting of both the coefficient and the exterior basis element.

```rust
pub fn basis_iter(&self) -> impl Iterator<Item = (f64, ExteriorTerm)> + '_ {
  let dim = self.dim;
  let grade = self.grade;
  self
    .coeffs
    .iter()
    .copied()
    .enumerate()
    .map(move |(i, coeff)| {
      let basis = ExteriorTerm::from_lex_rank(dim, grade, i);
      (coeff, basis)
    })
}
```

We then implemented the addition and scalar multiplication of exterior elements by just applying the operation to the coefficients.

## 3.5   Exterior Product

The most obvious operation on a `ExteriorElement` is of course the exterior product [5]. For this we rely on the exterior product of two `ExteriorTerm`s, which is just a concatenation of the two multiindices.

```rust
impl ExteriorTerm {
  pub fn wedge(mut self, mut other: Self) -> Self {
    self.indices.append(&mut other.indices);
    self
  }
}
```

For the `ExteriorElement` we just iterate over the all combinations of basis expansion and canonicalize the wedges of the individual terms.

```rust
impl ExteriorElement {
  pub fn wedge(&self, other: &Self) -> Self {
    assert_eq!(self.dim, other.dim);
    let dim = self.dim;

    let new_grade = self.grade + other.grade;
    if new_grade > dim {
      return Self::zero(dim, 0);
    }

    let new_basis_size = binomial(dim, new_grade);
    let mut new_coeffs = na::DVector::zeros(new_basis_size);

    for (self_coeff, self_basis) in self.basis_iter() {
      for (other_coeff, other_basis) in other.basis_iter() {
        let self_basis = self_basis.clone();

        let coeff_prod = self_coeff * other_coeff;
        if self_basis == other_basis || coeff_prod == 0.0 {
          continue;
        }
        if let Some((sign, merged_basis)) = self_basis.wedge(other_basis).canonicalized() {
          let merged_basis = merged_basis.lex_rank();
          new_coeffs[merged_basis] += sign.as_f64() * coeff_prod;
        }
      }
    }

    Self::new(new_coeffs, dim, new_grade)
  }
}
```

And we also implement a big wedge operator, that takes an iterator of factors.

```rust
pub fn wedge_big(factors: impl IntoIterator<Item = Self>) -> Option<Self> {
  let mut factors = factors.into_iter();
  let first = factors.next()?;
  let prod = factors.fold(first, |acc, factor| acc.wedge(&factor));
  Some(prod)
}
```

# 3.6  Inner product on Exterior Elements

For the weak formulations of our PDEs [1] we rely on Hilbert spaces that require an $L^2$-inner product on differential forms. This is derived directly from the point-wise inner product on multiforms. Which itself is derived from the inner product on the tangent space, which comes from the Riemannian metric at the point [5].

This derivation from the inner product on the tangent space $g_p$ to the inner product on the exterior fiber $\bigwedge^k T_p^* M$, shall be computed.

In general given an inner product on the linear space $V$, we can derive an inner product on the exterior power $\bigwedge^k(V)$. The rule is the following: [4]

$$\langle e_I, e_J \rangle = \det \left[ \left\langle \mathrm{d}x_{I_i}, \mathrm{d}x_{J_j} \right\rangle \right]_{i,j}^k \tag{40}$$

Computationally we represent inner products as Gramian matrices on some basis. This means that we compute an extended Gramian matrix as the inner product on multielements from the Gramian matrix of single elements using the determinant.

```rust
/// Construct Gramian on lexicographically ordered standard k-element standard
/// basis from Gramian on single elements.
pub fn multi_gramian(single_gramian: &Gramian, grade: ExteriorGrade) -> Gramian {
  let dim = single_gramian.dim();
  let bases: Vec<_> = exterior_bases(dim, grade).collect();

  let mut multi_gramian = Matrix::zeros(bases.len(), bases.len());
  let mut multi_basis_mat = Matrix::zeros(grade, grade);

  for icomb in 0..bases.len() {
    let combi = &bases[icomb];
    for jcomb in icomb..bases.len() {
      let combj = &bases[jcomb];

      for iicomb in 0..grade {
        let combii = combi[iicomb];
        for jjcomb in 0..grade {
          let combjj = combj[jjcomb];
          multi_basis_mat[(iicomb, jjcomb)] = single_gramian[(combii, combjj)];
        }
      }
      let det = multi_basis_mat.determinant();
      multi_gramian[(icomb, jcomb)] = det;
      multi_gramian[(jcomb, icomb)] = det;
    }
  }
  Gramian::new_unchecked(multi_gramian)
}
```

We are already at the end of the implementation of the exterior algebra. There exist many operations that could be implemented as well, such as the Hodge star operator $\star$ [5], [4], based on an inner product. While the Hodge star is theoretically central, the chosen mixed finite element formulation allows us to construct the discrete system using only the exterior derivative and the inner product, bypassing the need to explicitly assemble the Hodge star operator itself.

# Chapter 4

# Discrete Differential Forms: Simplicial Cochains and Whitney Forms

Having discretized the smooth domain $\Omega$ into a simplicial complex $\mathcal{M}$, we now require a corresponding discretization for the differential forms $\Lambda(\Omega)$ defined upon it [10]. FEEC relies on discrete counterparts, called discrete differential forms (DDF) $\Lambda_h(\mathcal{M})$ defined on the simplicial complex $\mathcal{M}$, that faithfully represent the structure of the continuous forms [4], [3].

Central to FEEC are two fundamental, closely related concepts: [10]

- **Simplicial Cochains**: These are the primary discrete objects. A $k$-cochain assigns a real value (a degree of freedom) to each $k$-simplex of the mesh $\mathcal{M}$. They form the algebraic backbone, capturing the combinatorial topology and enabling discrete versions of operators like the exterior derivative [8].
- **Whitney Forms**: These constitute the lowest-order finite element space for differential forms [10]. Each Whitney $k$-form is a piecewise polynomial basis function associated with a $k$-simplex. They provide a means to reconstruct a field across the mesh by interpolating the cochain values and are essential for defining the integrals required in weak formulations [4].

The space of $k$-cochains and the space spanned by Whitney $k$-forms are isomorphic. The **projection** from continuous forms to cochains is realized by the **integration map**, while the **interpolation** from cochains to the Whitney space is achieved via the **Whitney map** [10].

This chapter explores the representation and manipulation of these discrete differential forms. We will focus on two key processes:

- **Discretization**: Projecting continuous differential forms onto the discrete setting, yielding degrees of freedom associated with the mesh simplices.
- **Reconstruction**: Interpolating these discrete degrees of freedom using basis functions to obtain a piecewise continuous representation over the mesh.

Understanding these discrete structures and the maps connecting them is crucial, as they provide the foundation for constructing the finite element spaces and structure-preserving discrete operators used throughout FEEC [4].

## 4.1 Simplicial Cochains

The discretization of a differential $k$-form is a simplicial $k$-cochain. A simplicial $k$-cochain $\omega$ is a real-valued function $\omega : \Delta_k(\mathcal{M}) \to \mathbb{R}$ on the $k$-skeleton $\Delta_k(\mathcal{M})$ of the mesh $\mathcal{M}$ [10].

Simplicial cochains arise naturally from the combinatorial structure of a simplicial complex, they are the duals of simplicial chains. Simplicial cochains are also the fundamental combinatorial object in **discrete exterior calculus** (DEC) [8].

One can represent this function on the simplices, using a list of real values that are ordered according to the global numbering of the simplices.

```
pub struct Cochain {
  pub coeffs: na::DVector<f64>,
```

```rust
    pub dim: Dim,
}
```

Simplicial cochains preserve the cohomology of the de Rham complex at a discrete level and therefore retain the key topological properties [3].

## 4.1.1 Discretization: Cochain-Projection via Integration

What is the interpretation of these cochain values? This question can be answered by looking at the discretization procedure of a continuous differential form to this discrete form.

The discretization of differential forms happens by projection onto the cochain space. This cochain-projection is the simple operation of integrating the given continuous differential $k$-form over every $k$-simplex of the mesh. This gives a real number for each $k$-simplex, which is exactly the cochain values. This projection is called the **integration map** [10]

$$I : \Lambda^k(\Omega) \to C^k(\mathcal{M}; \mathbb{R}) \tag{41}$$

$$I(\omega) = (\sigma \mapsto c_\sigma) \quad \text{where} \quad c_\sigma = \int_\sigma \omega \quad \forall \sigma \in \Delta_k(\mathcal{M}) \tag{42}$$

The integral of a differential $k$-form over a $k$-simplex is defined using the pullback to the reference $k$-simplex.

$$
\begin{aligned}
\int_\sigma \omega &= \int_{\hat{\sigma}} \varphi^* \omega \\
&= \int_{\hat{\sigma}} \omega_{\varphi(\lambda^1,\dots,\lambda^k)}(\varphi_* \hat{e}_1, \dots, \varphi_* \hat{e}_k) \mathrm{d}\lambda^1 \dots \mathrm{d}\lambda^k \\
&= \int_{\hat{\sigma}} \omega_{\varphi(\lambda^1,\dots,\lambda^k)}(\boldsymbol{e}_1, \dots, \boldsymbol{e}_n) \mathrm{d}\lambda^1 \dots \mathrm{d}\lambda^k
\end{aligned}
\tag{43}
$$

The last expression is a traditional pre-differential geometry integral over a subset $\hat{\sigma} \in \mathbb{R}^k$. No exterior calculus required. We approximate this ordinary integral using barycentric quadrature.

$$\int_\sigma \omega \approx |\hat{\sigma}| \; \omega_{\varphi(\boldsymbol{m}_\sigma)}(\boldsymbol{e}_1, \dots, \boldsymbol{e}_n) \tag{44}$$

Here the $e_i$ are the spanning vectors of the simplex $\sigma$. We work here with coordinate representation of differential froms, that rely on an embedding. Computationally we will have functors that can be evaluated at coordinates to get the multiform.

And the implementation just looks like this:

```rust
pub fn integrate_form_simplex(form: &impl DifferentialMultiForm, simplex: &SimplexCoords) -> f64 {
  let multivector = simplex.spanning_multivector();
  let f = |coord: CoordRef| {
    form
      .at_point(simplex.local2global(coord).as_view())
      .apply_form_on_multivector(&multivector)
  };
  let std_simp = SimplexCoords::standard(simplex.dim_intrinsic());
  barycentric_quadrature(&f, &std_simp)
}
```

And for the full cochain-projection, we just repeat this integration for each simplex in the $k$-skeleton. The implementation the looks like this.

```rust
pub fn cochain_projection(
  form: &impl DifferentialMultiForm,
  topology: &Complex,
  coords: &MeshCoords,
) -> Cochain {
  let cochain = topology
    .skeleton(form.grade())
    .handle_iter()
```

```
    .map(|simp| SimplexCoords::from_simplex_and_coords(&simp, coords))
    .map(|simp| integrate_form_simplex(form, &simp))
    .collect::<Vec<_>>()
    .into();
  Cochain::new(form.grade(), cochain)
}
```

## 4.1.2 Discrete Exterior Derivative via Stokes' Theorem

In exterior calculus, the exterior derivative is a fundamental operator that generalizes the standard derivatives (gradient, curl, divergence) to differential forms [5].

$$\mathrm{d} : \Lambda^k(\Omega) \to \Lambda^{k+1}(\Omega) \tag{45}$$

For our discrete setting, we require a discrete counterpart that acts on cochains. A so called **discrete exterior derivative** $\mathrm{d}_h$, which can be derived through the lens of cochain calculus.

$$\mathrm{d}_h : C^k(\mathcal{M}) \to C^{k+1}(\mathcal{M}) \tag{46}$$

A crucial property of the continuous exterior derivative is captured by **Stokes' Theorem on chains**. For a differential form $\omega$ and a chain $c$, this theorem relates the exterior derivative with the boundary operator [5], [6]:

$$\int_c \mathrm{d}\omega = \int_{\partial c} \omega \tag{47}$$

This relationship is particularly insightful when viewed through the framework of dual pairings. If we define a natural pairing $\langle \cdot, \cdot \rangle$ between differential forms and chains as integration over the chain:

$$\langle \omega, c \rangle := \int_c \omega \tag{48}$$

Now Stokes' Theorem can be expressed in a more abstract form:

$$\langle \mathrm{d}\omega, c \rangle = \langle \omega, \partial c \rangle \tag{49}$$

This equation reveals that, with respect to this dual pairing, the exterior derivative operator $\mathrm{d}$ is the adjoint of the boundary operator $\partial$.

$$\mathrm{d} = \partial^* \tag{50}$$

This adjoint relationship provides the direct motivation for defining the discrete exterior derivative. The discrete exterior derivative, often called the **coboundary operator**, is thus defined as the adjoint of the boundary operator. By definition, this discrete operator preserves the structure of Stokes' Theorem at the discrete level [8].

From a computational perspective, the boundary operator $\partial^k$ mapping $k$-chains to $(k-1)$-chains can be represented as a signed incidence matrix between the simplices in the $k$-skeleton and $(k-1)$-skeleton [8]. The adjoint property then translates directly to the discrete exterior derivative $\mathrm{d}^k$ (mapping $k$-cochains to $(k+1)$-cochains) being the transpose of the boundary operator $\partial_{k+1}$:

$$\mathbf{d}^k = \mathbf{D}_{k+1}^{\top}$$
$$\mathbf{d}^k \in \{-1, 0, +1\}^{N_{k+1} \times N_k} \tag{51}$$

This definition highlights a key feature of the discrete exterior derivative: it is a purely topological operator. Its definition and computation depend only on the combinatorial structure of the simplicial complex (captured by the incidence relationships in the boundary operator), and not on the geometry or metric of the underlying manifold [8].

In our implementation, we represent the discrete exterior derivative as a sparse matrix. Using an extension trait, we provide a method to compute this matrix for a given grade, leveraging the already implemented boundary operator:

```
pub trait ManifoldComplexExt { ... }
impl ManifoldComplexExt for Complex {
  /// $dif^k: cal(W) Lambda^k -> cal(W) Lambda^(k+1)$
  fn exterior_derivative_operator(&self, grade: ExteriorGrade) -> SparseMatrix {
    self.boundary_operator(grade + 1).transpose()
```

}
    }

## 4.2 Whitney Forms

**Whitney forms** are the **finite element differential forms** that correspond to cochains. They are the piecewise-linear (over the cells) differential forms defined over the simplicial manifold [10].

The Whitney space $\mathcal{W}\Lambda^k(\mathcal{M})$ is the space of all Whitney forms over our mesh $\mathcal{M}$

The **Whitney subcomplex** is a subcomplex of the continuous de Rham complex of differential forms. [4]

$$0 \to \mathcal{W}\Lambda^0(\mathcal{M}) \xrightarrow{\mathrm{d}} \cdots \xrightarrow{\mathrm{d}} \mathcal{W}\Lambda^n(\mathcal{M}) \to 0 \tag{52}$$

### 4.2.1 Whitney Basis

There is a special basis for the space of Whitney forms, called the **Whitney basis** [10], [3]. Just like there is a cochain value for each $k$-simplex, there is a Whitney basis function for each $k$-simplex. They have their DOF on this $k$-simplex.

$$\mathcal{W}\Lambda^k(\mathcal{M}) = \mathrm{span}\left\{\varphi_\sigma : \sigma \in \Delta_k(\mathcal{M})\right\} \tag{53}$$

Let's take a look at the local shape functions (LSF).

The local Whitney form $\varphi_\sigma|_K = \lambda_\sigma = \lambda_{i_0...i_k}$ associated with the DOF simplex $\sigma = [i_0...i_k] \subseteq K$ on the cell $K = [j_0...j_n]$ is defined using the barycentric coordinate functions $\lambda_{j_s} : K \to \mathbb{R}$ of the cell.

$$\lambda_{i_0...i_k} = k! \sum_{l=0}^{k} (-1)^l \lambda_{i_l}\left(\mathrm{d}\lambda_{i_0} \wedge \cdots \wedge \widehat{\mathrm{d}\lambda_{i_l}} \wedge \cdots \wedge \mathrm{d}\lambda_{i_k}\right) \tag{54}$$

To get some intuition for the kind of fields this produces, let's look at some visualizations. We do this specifically for 1-forms, since we can intuitively visualize these using vector field proxies.

We have the following formula for local shape functions of Whitney 1-forms.

$$\lambda_{ij} = \lambda_i \mathrm{d}\lambda_j - \lambda_j \mathrm{d}\lambda_i \tag{55}$$

For the reference 2-simplex, we get the following Whitney basis 1-forms.

$$\begin{aligned} \lambda_{01} &= (1-y)\mathrm{d}x + x\mathrm{d}y \\ \lambda_{02} &= y\mathrm{d}x + (1-x)\mathrm{d}y \\ \lambda_{12} &= -y\mathrm{d}x + x\mathrm{d}y \end{aligned} \tag{56}$$

Visualized on the reference triangle they look like:

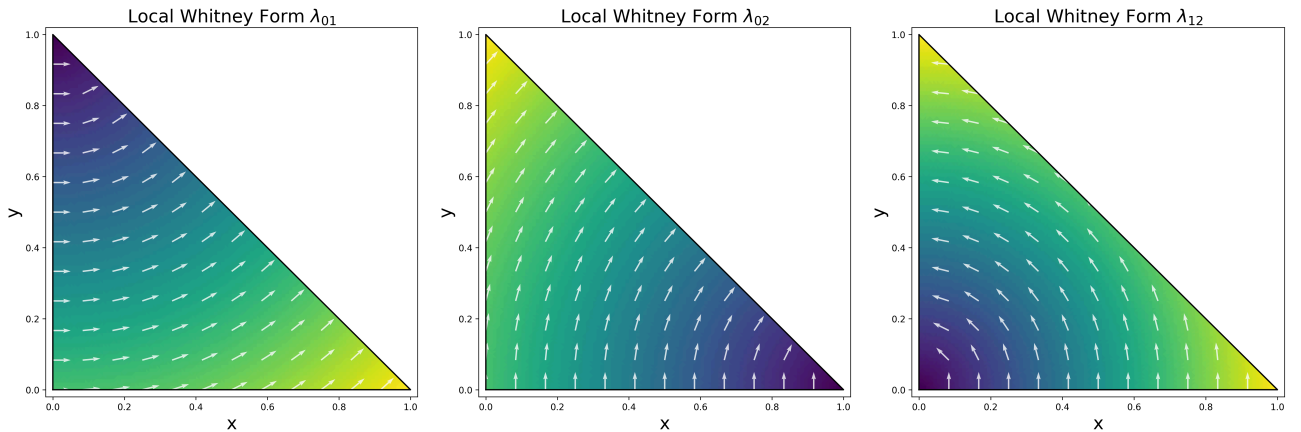

Figure 1: Vector proxies of Reference Local Shape Functions $\lambda_{01}, \lambda_{02}, \lambda_{12} \in \mathcal{W}\Lambda^1\left(\Delta_2^{\mathrm{ref}}\right)$.

The global shape functions are obtained by combing the various local shape functions that are associated with the same DOF simplex. In general the GSF are discontinuous over cell boundaries.

We visualize here the 3 GSF on a equilateral triangle mesh associated with the edges of the middle triangle.

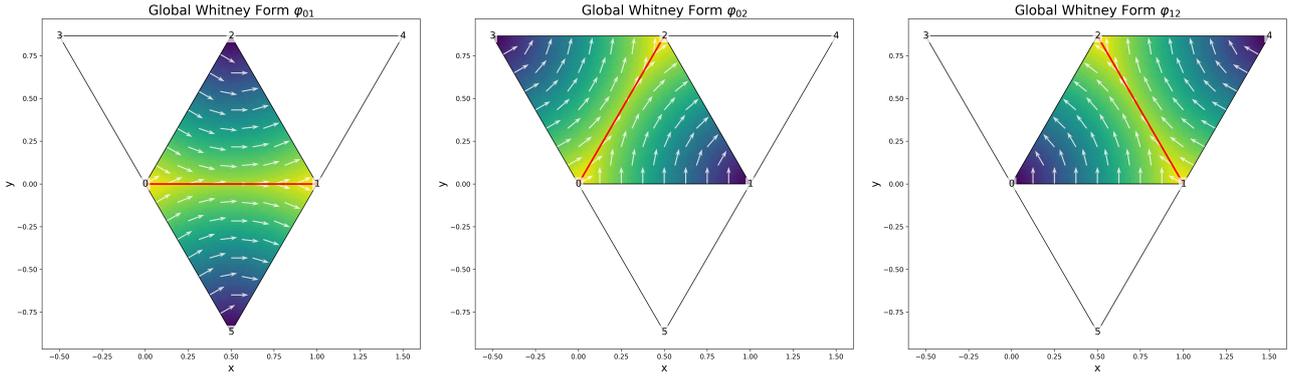

Figure 2: Vector proxies of Global Shape Functions $\varphi_{01}, \varphi_{02}, \varphi_{12} \in \mathcal{W}\Lambda^1(\mathcal{M})$ on equilateral triangle mesh $\mathcal{M}$.

Via linear combination of these GSF we can obtain any Whitney form $u = \sum_{\sigma \in \Delta_k(\mathcal{M})} c^\sigma \varphi_\sigma$, where $c \in C^k(\mathcal{M})$ is a simplicial cochain providing the coefficients $c^i$.

We can for instance construct the following constant, purely divergent, and purely rotational vector fields.

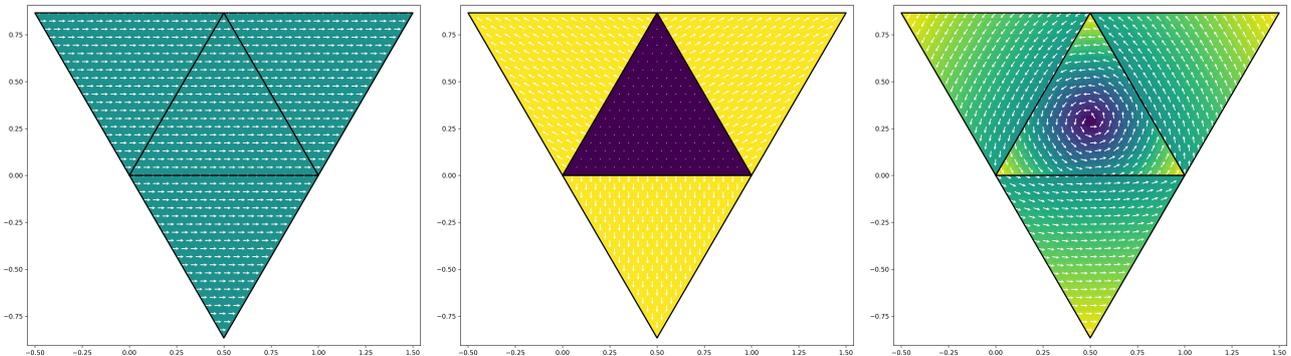

Figure 3: Vector proxies of some example Whitney forms on equilateral triangle mesh $\mathcal{M}$.

We now implement some functionality to work with Whitney forms. The most important of which is the point evaluation of Whitney basis forms on a cell associated with one of it's subsimplices. Which means implementing Equation 54.

We do this implementation based on a coordinate simplex. But this doesn't mean that this relies on an embedding. It is multi-purpose! If we just choose the reference simplex as the coordinate simplex, then we just compute the reference Whitney basis forms. If we choose a coordinate simplex from an embedding, we compute the representation of a Whitney form in this embedding.

We define a struct that stores the coordinate of the cell we are working on as well as the DOF subsimplex (in local indices) that the Whitney basis function is associated with.

```rust
#[derive(Debug, Clone)]
pub struct WhitneyLsf {
  cell_coords: SimplexCoords,
  dof_simp: Simplex,
}
impl WhitneyLsf {
  pub fn from_coords(cell_coords: SimplexCoords, dof_simp: Simplex) -> Self {
    Self {
      cell_coords,
      dof_simp,
    }
  }
  pub fn standard(cell_dim: Dim, dof_simp: Simplex) -> Self {
    Self::from_coords(SimplexCoords::standard(cell_dim), dof_simp)
  }
  pub fn grade(&self) -> ExteriorGrade { self.dof_simp.dim() }
```

The defining formula of a Whitney basis form, relies on the barycentric coordinate functions of the cell, but only those that are on the vertices of the DOF simplex. We write a function to get an iterator on exactly these.

```rust
/// The difbarys of the vertices of the DOF simplex.
pub fn difbarys(&self) -> impl Iterator<Item = MultiForm> + use<'_> {
  self
    .cell_coords
    .difbarys_ext()
    .into_iter()
    .enumerate()
    .filter_map(|(ibary, difbary)| self.dof_simp.contains(ibary).then_some(difbary))
}
```

We can observe that the big wedge terms $d\lambda_{i_0} \wedge \cdots \wedge \widehat{d\lambda_{i_l}} \wedge \cdots \wedge d\lambda_{i_k}$ in Equation 54 are constant multiforms. We write a function that computes these constants.

```rust
/// dλ_i_0 ∧⋯∧^omit(dλ_i_iwedge) ∧⋯∧ dλ_i_dim
pub fn wedge_term(&self, iterm: usize) -> MultiForm {
  let dim_cell = self.cell_coords.dim_intrinsic();
  let wedge = self
    .difbarys()
    .enumerate()
    // leave off i'th difbary
    .filter_map(|(pos, difbary)| (pos != iterm).then_some(difbary));
  MultiForm::wedge_big(wedge).unwrap_or(MultiForm::one(dim_cell))
}
pub fn wedge_terms(&self) -> impl ExactSizeIterator<Item = MultiForm> + use<'_> {
  (0..self.dof_simp.nvertices()).map(move |iwedge| self.wedge_term(iwedge))
}
```

Now we can implement the `ExteriorField` trait for our `WhitneyLsf` struct to make it point evaluable, based on a coordinate.

```rust
impl ExteriorField for WhitneyLsf {
  fn dim_ambient(&self) -> exterior::Dim { self.cell_coords.dim_ambient() }
  fn dim_intrinsic(&self) -> exterior::Dim { self.cell_coords.dim_intrinsic() }
  fn grade(&self) -> ExteriorGrade { self.grade() }
  fn at_point<'a>(&self, coord: impl Into<CoordRef<'a>>) -> MultiForm {
    let barys = self.cell_coords.global2bary(coord);
    assert!(is_bary_inside(&barys), "Point is outside cell.");

    let dim = self.dim_intrinsic();
    let grade = self.grade();
    let mut form = MultiForm::zero(dim, grade);
    for (iterm, &vertex) in self.dof_simp.vertices.iter().enumerate() {
      let sign = Sign::from_parity(iterm);
      let wedge = self.wedge_term(iterm);

      let bary = barys[vertex];
      form += sign.as_f64() * bary * wedge;
    }
    (factorial(grade) as f64) * form
  }
}
```

We can implement one more functionality.

Since the Whitney form is a linear differential form over the cell, it's exterior derivative must be a constant multiform. We can easily derive it's value.

$$
\begin{aligned}
\mathrm{d}\lambda_{i_0\cdots i_k} &= k! \sum_{l=0}^{k} (-1)^l \mathrm{d}\lambda_{i_l} \wedge \left( \mathrm{d}\lambda_{i_0} \wedge \cdots \wedge \widehat{\mathrm{d}\lambda}_{i_l} \wedge \cdots \wedge \mathrm{d}\lambda_{i_k} \right) \\
&= k! \sum_{l=0}^{k} (-1)^l (-1)^l \left( \mathrm{d}\lambda_{i_0} \wedge \cdots \wedge \mathrm{d}\lambda_{i_l} \wedge \cdots \wedge \mathrm{d}\lambda_{i_k} \right) \\
&= (k+1)! \mathrm{d}\lambda_{i_0} \wedge \cdots \wedge \mathrm{d}\lambda_{i_k}
\end{aligned}
\tag{57}
$$

This is the corresponding implementation.

```rust
pub fn dif(&self) -> MultiForm {
  let dim = self.cell_coords.dim_intrinsic();
  let grade = self.grade();
  if grade == dim {
    return MultiForm::zero(dim, grade + 1);
  }
  factorialf(grade + 1) * MultiForm::wedge_big(self.difbarys()).unwrap()
}
```

The defining property of the Whitney basis is a from pointwise to integral generalized Lagrange basis property [10]:
For any two $k$-simplices $\sigma, \tau \in \Delta_k(\mathcal{M})$, we have

$$
\int_\sigma \lambda_\tau = \begin{cases} +1 & \text{if } \sigma = +\tau \\ -1 & \text{if } \sigma = -\tau \\ 0 & \text{if } \sigma \neq \pm\tau \end{cases}
\tag{58}
$$

We can write a test that verifies our implementation by checking this property.

```rust
#[test]
fn whitney_basis_property() {
  for dim in 0..=4 {
    let topology = Complex::standard(dim);
    let coords = MeshCoords::standard(dim);

    for grade in 0..=dim {
      for dof_simp in topology.skeleton(grade).handle_iter() {
        let whitney_form = WhitneyLsf::standard(dim, (*dof_simp).clone());

        for other_simp in topology.skeleton(grade).handle_iter() {
          let are_same_simp = dof_simp == other_simp;
          let other_simplex = other_simp.coord_simplex(&coords);
          let discret = integrate_form_simplex(&whitney_form, &other_simplex);
          let expected = are_same_simp as u8 as f64;
          let diff = (discret - expected).abs();
          const TOL: f64 = 10e-9;
          let equal = diff <= TOL;
          assert!(equal, "for: computed={discret} expected={expected}");
          if other_simplex.nvertices() >= 2 {
            let other_simplex_rev = other_simplex.clone().flipped_orientation();
            let discret_rev = integrate_form_simplex(&whitney_form, &other_simplex_rev);
            let expected_rev = Sign::Neg.as_f64() * are_same_simp as usize as f64;
            let diff_rev = (discret_rev - expected_rev).abs();
            let equal_rev = diff_rev <= TOL;
            assert!(
              equal_rev,
              "rev: computed={discret_rev} expected={expected_rev}"
            );
          }
        }
      }
    }
  }
}
```

### 4.2.2 Reconstruction: Whitney-Interpolation via the Whitney map

There is a one-to-one correspondence between $k$-cochains and Whitney $k$-forms [10]. We have already seen how to obtain a cochain from a continuous differential form via cochain-projection. But what about the other way? How can we reconstruct a continuous differential form from a cochain?

This reconstruction is achieved by the so called **Whitney map** [10].

$$W : C^k(\mathcal{M}; \mathbb{R}) \to \mathcal{W}^k(\mathcal{M}) \subseteq \Lambda^k(\Omega)$$
$$W(c) = \sum_{\sigma \in \mathcal{M}_k} c(\sigma)\omega_\sigma \tag{59}$$

It can be seen as a generalized interpolation operator. Instead of pointwise interpolation, we have interpolation in an integral sense. It takes cochains to differential forms that have the cochains values as integral values [10].

The isomorphism between Whitney forms and cochains can now be constructed

$$W \circ I = \mathrm{id}_{\mathcal{W}^k(\mathcal{M})}$$
$$I \circ W = \mathrm{id}_{C^k(\mathcal{M}; \mathbb{R})} \tag{60}$$

On our implementation side, we introduce a routine that takes a cochain and let's us evaluate the corresponding Whitney form at any point on the simplicial manifold.

```
pub fn whitney_form_eval<'a>(
  coord: impl Into<CoordRef<'a>>,
  cochain: &Cochain,
  mesh_cell: SimplexHandle,
  mesh_coords: &MeshCoords,
) -> MultiForm {
  let coord = coord.into();

  let cell_coords = mesh_cell.coord_simplex(mesh_coords);

  let dim_intrinsic = mesh_cell.dim();
  let grade = cochain.dim;

  let mut value = MultiForm::zero(dim_intrinsic, grade);
  for dof_simp in mesh_cell.mesh_subsimps(grade) {
    let local_dof_simp = dof_simp.relative_to(&mesh_cell);

    let lsf = WhitneyLsf::from_coords(cell_coords.clone(), local_dof_simp);
    let lsf_value = lsf.at_point(coord);
    let dof_value = cochain[dof_simp];
    value += dof_value * lsf_value;
  }
  value
}
```

We shall provide the cell on which the point lives, to avoid searching over all the cells for the one containing the point.

## 4.3 Higher-Order Discrete Differential Forms

The theoretical construction of finite element differential forms exist for all polynomial degrees [4], [2]. We don't support them in this implementation, but this a very obvious possible future extension to this implementation. One just needs to keep in mind that then higher-order manifold approximations are also needed to not incur any non-admissible geometric variational crimes [11].

# Chapter 5

# Finite Element Methods for Differential Forms

We have now arrived at the chapter discussing the actual `formoniq` crate.

Here we will derive and implement the formulas for computing the element matrices of the various weak differential operators in Finite Element Exterior Calculus (FEEC) [3], [4]. Furthermore we implement the assembly algorithm that will give us the final Galerkin matrices [1].

## 5.1 Variational Formulation & Element Matrix Computation

There are only 4 types of variational operators that are relevant to the mixed weak formulation of the Hodge-Laplace operator [4]. All of them are based on the $L^2$ inner product on Whitney forms [10], [3].

Above all is the mass bilinear form, which is the $L^2$ inner product [4].

$$m^k(u,v) = \langle u, v \rangle_{L^2\Lambda^k(\Omega)} \quad u \in L^2\Lambda^k, v \in L^2\Lambda^k \tag{61}$$

The next bilinear form involves the exterior derivative.

$$d^k(u,v) = \langle \mathrm{d}u, v \rangle_{L^2\Lambda^k(\Omega)} \quad u \in H\Lambda^{k-1}, v \in L^2\Lambda^k \tag{62}$$

The bilinear form involving the codifferential is also relevant [3].

$$c(u,v) = \langle \delta u, v \rangle_{L^2\Lambda^k(\Omega)} \tag{63}$$

Using the adjoint property of the codifferential relative to the exterior derivative under the $L^2$ inner product [4], it can be rewritten using the exterior derivative applied to the test function.

$$c^k(u,v) = \langle u, \mathrm{d}v \rangle_{L^2\Lambda^k(\Omega)} \quad u \in L^2\Lambda^k, v \in H\Lambda^{k-1} \tag{64}$$

Lastly there is the bilinear form analogous to the scalar Laplacian, involving exterior derivatives on both arguments. It corresponds to the $\delta\mathrm{d}$ part of the Hodge-Laplacian [4].

$$l^k(u,v) = \langle \mathrm{d}u, \mathrm{d}v \rangle_{L^2\Lambda^{k+1}(\Omega)} \quad u \in H\Lambda^k, v \in H\Lambda^k \tag{65}$$

After Galerkin discretization [1], by means of the Whitney finite element space [10] with the Whitney basis $\{\varphi_i^k\}$, we arrive at the following Galerkin matrices for our four weak operators.

$$\begin{aligned}
\mathbf{M}^k &= \left[ \langle \varphi_i^k, \varphi_j^k \rangle \right]_{ij} \\
\mathbf{D}^k &= \left[ \langle \varphi_i^k, \mathrm{d}\varphi_j^{k-1} \rangle \right]_{ij} \\
\mathbf{C}^k &= \left[ \langle \mathrm{d}\varphi_i^{k-1}, \varphi_j^k \rangle \right]_{ij} \\
\mathbf{L}^k &= \left[ \langle \mathrm{d}\varphi_i^k, \mathrm{d}\varphi_j^k \rangle \right]_{ij}
\end{aligned} \tag{66}$$

We can rewrite the 3 operators that involve the exterior derivative using the mass matrix and the discrete exterior derivative (coboundary/incidence matrix) [3], [8].

$$\mathbf{D}^k = \mathbf{M}^k \mathbf{d}^{k-1}$$
$$\mathbf{C}^k = \left(\mathbf{d}^{k-1}\right)^\top \mathbf{M}^k \tag{67}$$
$$\mathbf{L}^k = \left(\mathbf{d}^k\right)^\top \mathbf{M}^{k+1} \mathbf{d}^k$$

As usual in a FEM library [1], we define element matrix providers, that compute the element matrices on each cell of mesh and later on assemble the full Galerkin matrices.

We first define a element matrix provider trait

```rust
pub type ElMat = Matrix;
pub trait ElMatProvider: Sync {
  fn row_grade(&self) -> ExteriorGrade;
  fn col_grade(&self) -> ExteriorGrade;
  fn eval(&self, geometry: &SimplexLengths) -> ElMat;
}
```

The eval method provides us with the element matrix on a cell, given it's geometry. But we also need to know the exterior grade of the Whitney forms that correspond to the rows and columns. This information will be used by the assembly routine.

We will now implement the 3 operators involving exterior derivatives based on the element matrix provider of the mass bilinear form.

The local exterior derivative only depends on the local topology, which is the same for any simplex of the same dimension. So we use a global variable that stores the transposed incidence matrix (coboundary operator) for any $k$-skeleton of a $n$-complex [8].

```rust
pub struct DifElmat {
  mass: HodgeMassElmat,
  dif: Matrix,
}
impl DifElmat {
  pub fn new(dim: Dim, grade: ExteriorGrade) -> Self {
    let mass = HodgeMassElmat::new(dim, grade);
    let dif = Complex::standard(dim).exterior_derivative_operator(grade - 1);
    let dif = Matrix::from(&dif);
    Self { mass, dif }
  }
}
impl ElMatProvider for DifElmat {
  fn row_grade(&self) -> ExteriorGrade { self.mass.grade }
  fn col_grade(&self) -> ExteriorGrade { self.mass.grade - 1 }
  fn eval(&self, geometry: &SimplexLengths) -> Matrix {
    let mass = self.mass.eval(geometry);
    mass * &self.dif
  }
}

pub struct CodifElmat {
  mass: HodgeMassElmat,
  codif: Matrix,
}
impl CodifElmat {
  pub fn new(dim: Dim, grade: ExteriorGrade) -> Self {
    let mass = HodgeMassElmat::new(dim, grade);
    let dif = Complex::standard(dim).exterior_derivative_operator(grade - 1);
    let dif = Matrix::from(&dif);
    let codif = dif.transpose();
    Self { mass, codif }
  }
}
```

```rust
impl ElMatProvider for CodifElmat {
  fn row_grade(&self) -> ExteriorGrade { self.mass.grade - 1 }
  fn col_grade(&self) -> ExteriorGrade { self.mass.grade }
  fn eval(&self, geometry: &SimplexLengths) -> Matrix {
    let mass = self.mass.eval(geometry);
    &self.codif * mass
  }
}

pub struct CodifDifElmat {
  mass: HodgeMassElmat,
  dif: Matrix,
  codif: Matrix,
}
impl CodifDifElmat {
  pub fn new(dim: Dim, grade: ExteriorGrade) -> Self {
    let mass = HodgeMassElmat::new(dim, grade + 1);
    let dif = Complex::standard(dim).exterior_derivative_operator(grade);
    let dif = Matrix::from(&dif);
    let codif = dif.transpose();
    Self { mass, dif, codif }
  }
}
impl ElMatProvider for CodifDifElmat {
  fn row_grade(&self) -> ExteriorGrade { self.mass.grade - 1 }
  fn col_grade(&self) -> ExteriorGrade { self.mass.grade - 1 }
  fn eval(&self, geometry: &SimplexLengths) -> Matrix {
    let mass = self.mass.eval(geometry);
    &self.codif * mass * &self.dif
  }
}
```

### 5.1.1 Mass bilinear form

Now we need to implement the element matrix provider to the mass bilinear form. Here is where the geometry of the domain comes in, through the inner product, which depends on the Riemannian metric [5].

One could also understand the mass bilinear form as a weak Hodge star operator [4].

$$\mathbf{M}_{ij} = \int_\Omega \varphi_j \wedge \star \varphi_i = \langle \varphi_j, \varphi_i \rangle_{L^2 \Lambda^k(\Omega)} \tag{68}$$

We will not compute this using the Hodge star operator, but instead directly using the inner product.

We already have an inner product on constant multiforms. We now need to extend it to an $L^2$ inner product on Whitney forms [3]. This can be done by inserting the definition of a Whitney form (in terms of barycentric coordinate functions) into the inner product.

$$
\begin{aligned}
\langle \lambda_{i_0 \dots i_k}, \lambda_{j_0 \dots j_k} \rangle_{L^2 \Lambda^k(\Omega)} &= k!^2 \sum_{l=0}^{k} \sum_{m=0}^{k} (-)^{l+m} \left\langle \begin{matrix} \lambda_{i_l} \left( \mathrm{d}\lambda_{i_0} \wedge \cdots \wedge \widehat{\mathrm{d}\lambda}_{i_l} \wedge \cdots \wedge \mathrm{d}\lambda_{i_k} \right) \\ \lambda_{j_m} \left( \mathrm{d}\lambda_{j_0} \wedge \cdots \wedge \widehat{\mathrm{d}\lambda}_{j_m} \wedge \cdots \wedge \mathrm{d}\lambda_{j_k} \right) \end{matrix} \right\rangle_{L^2 \Lambda^k(\Omega)} \\
&= k!^2 \sum_{l,m} (-)^{l+m} \left\langle \begin{matrix} \mathrm{d}\lambda_{i_0} \wedge \cdots \wedge \widehat{\mathrm{d}\lambda}_{i_l} \wedge \cdots \wedge \mathrm{d}\lambda_{i_k} \\ \mathrm{d}\lambda_{j_0} \wedge \cdots \wedge \widehat{\mathrm{d}\lambda}_{j_m} \wedge \cdots \wedge \mathrm{d}\lambda_{j_k} \end{matrix} \right\rangle_{\Lambda^k} \int_K \lambda_{i_l} \lambda_{j_m} \mathrm{vol}_g
\end{aligned}
\tag{69}
$$

We can now make use of the fact that the exterior derivative of the barycentric coordinate functions are constant [1]. This makes the big wedge terms also constant. We can therefore pull them out of the $L^2$ integral and now it's just an inner product on constant multiforms. What remains in the in the integral is the product of two barycentric coordinate functions.

Using this we can now implement the element matrix provider to the mass bilinear form in Rust.

```rust
pub struct HodgeMassElmat {
  dim: Dim,
  grade: ExteriorGrade,
```

```rust
  simplices: Vec<Simplex>,
  wedge_terms: Vec<ExteriorElementList>,
}
impl HodgeMassElmat {
  pub fn new(dim: Dim, grade: ExteriorGrade) -> Self {
    let simplices: Vec<_> = standard_subsimps(dim, grade).collect();
    let wedge_terms: Vec<ExteriorElementList> = simplices
      .iter()
      .cloned()
      .map(|simp| WhitneyLsf::standard(dim, simp).wedge_terms().collect())
      .collect();

    Self { dim, grade, simplices, wedge_terms }
  }
}
impl ElMatProvider for HodgeMassElmat {
  fn row_grade(&self) -> ExteriorGrade { self.grade }
  fn col_grade(&self) -> ExteriorGrade { self.grade }

  fn eval(&self, geometry: &SimplexLengths) -> Matrix {
    assert_eq!(self.dim, geometry.dim());

    let scalar_mass = ScalarMassElmat.eval(geometry);

    let mut elmat = Matrix::zeros(self.simplices.len(), self.simplices.len());
    for (i, asimp) in self.simplices.iter().enumerate() {
      for (j, bsimp) in self.simplices.iter().enumerate() {
        let wedge_terms_a = &self.wedge_terms[i];
        let wedge_terms_b = &self.wedge_terms[j];
        let wedge_inners = multi_gramian(&geometry.to_metric_tensor().inverse(), self.grade)
          .inner_mat(wedge_terms_a.coeffs(), wedge_terms_b.coeffs());

        let nvertices = self.grade + 1;
        let mut sum = 0.0;
        for avertex in 0..nvertices {
          for bvertex in 0..nvertices {
            let sign = Sign::from_parity(avertex + bvertex);
            let inner = wedge_inners[(avertex, bvertex)];
            sum += sign.as_f64() * inner * scalar_mass[(asimp[avertex], bsimp[bvertex])];
          }
        }
        elmat[(i, j)] = sum;
      }
    }
    factorial(self.grade).pow(2) as f64 * elmat
  }
}
```

Now we are just missing an element matrix provider for the scalar mass bilinear form. Luckily there exists a closed form solution, for this integral, which only depends on the volume of the cell [1].

$$\int_K \lambda_i \lambda_j \ \mathrm{vol} = \frac{|K|}{(n+2)(n+1)}\big(1+\delta_{ij}\big) \tag{70}$$

derived from this more general integral formula for powers of barycentric coordinate functions [1]

$$\int_K \lambda_0^{\alpha_0} \cdots \lambda_n^{\alpha_n} \ \mathrm{vol} = n!|K|\frac{\alpha_0! \cdots \alpha_n!}{(\alpha_0 + \cdots + \alpha_n + n)!} \tag{71}$$

```rust
pub struct ScalarMassElmat;
impl ElMatProvider for ScalarMassElmat {
  fn row_grade(&self) -> ExteriorGrade { 0 }
  fn col_grade(&self) -> ExteriorGrade { 0 }
  fn eval(&self, geometry: &SimplexGeometry) -> ElMat {
    let ndofs = geometry.nvertices();
```

```
    let dim = geometry.dim();
    let v = geometry.vol() / ((dim + 1) * (dim + 2)) as f64;
    let mut elmat = na::DMatrix::from_element(ndofs, ndofs, v);
    elmat.fill_diagonal(2.0 * v);
    elmat
  }
}
```

## 5.2  Source Terms & Quadrature

Next to bilinear forms and element matrix providers, there are also linear forms and corresponding element vector providers.

We define a trait for element vector providers, just like we did for the element matrices.

```
pub type ElVec = Vector;
pub trait ElVecProvider: Sync {
  fn grade(&self) -> ExteriorGrade;
  fn eval(&self, geometry: &SimplexLengths, topology: &Simplex) -> ElVec;
}
```

Here in contrast to the element matrix provider, we also pass the topology as an argument, to the evaluation method. This might come as a surprise, as this element vector should be purely local and therefore not work with any global topological information. However the element matrix provider for the source term will need global coordinates, but since we don't want to include coordinates in our `ElVecProvider` trait, we need to supply information on which topological simplex we are on, in order to reconstruct the corresponding coordinates.

We only have a single implementor of the `ElVecProvider` trait, which is the element vector that corresponds to the right-hand-side source term, as it appears for instance in the Hodge-Laplace source problem.

$$\boldsymbol{b}_K = \left[\left\langle f, \varphi_j^k \right\rangle_{L^2\Lambda^k(K)}\right]_{i=1}^{N_k} \tag{72}$$

It takes as input an `ExteriorField` that corresponds to $f$, represented as a functor that takes a global coordinate $\boldsymbol{x} \in \mathbb{R}^N$ (based on an embedding) and produces the multiform $f_{\boldsymbol{x}}$ at that position. It constitutes the coordinate-based representation of the arbitrary continuous differential form $f$.

Since only point-evaluation is available to us, we need to rely on quadrature [1] to compute the element vector corresponding to the source term. For this we pullback the `ExteriorField` representing f to the reference cell using `precompose_form` and compute the pointwise inner product.

$$\left\langle f, \varphi_j^k \right\rangle_{L^2\Lambda^k} = \int_K \left\langle f(\boldsymbol{x}), \varphi_j^k(\boldsymbol{x}) \right\rangle_{\Lambda^k} \mathrm{vol}_g \approx |K| \sum_l w_l \left\langle \varphi_K^\star f(\boldsymbol{q}_l), \varphi_j^k(\boldsymbol{q}_l) \right\rangle_{\Lambda^k} \tag{73}$$

The choice of quadrature is free, but simple barycentric quadrature suffices, to get an admissible variational crime for 1st order FE. [1]

This is what the full implementation of this element vector provider then looks like.

```
pub struct SourceElVec<'a, F>
where
  F: ExteriorField,
{
  source: &'a F,
  mesh_coords: &'a MeshCoords,
  qr: SimplexQuadRule,
}
impl<'a, F> SourceElVec<'a, F>
where
  F: ExteriorField,
{
  pub fn new(source: &'a F, mesh_coords: &'a MeshCoords, qr: Option<SimplexQuadRule>) -> Self {
    let qr = qr.unwrap_or(SimplexQuadRule::barycentric(source.dim_intrinsic()));
    Self { source, mesh_coords, qr }
```

```rust
  }
}
impl<F> ElVecProvider for SourceElVec<'_, F>
where
  F: Sync + ExteriorField,
{
  fn grade(&self) -> ExteriorGrade {
    self.source.grade()
  }
  fn eval(&self, geometry: &SimplexLengths, topology: &Simplex) -> ElVec {
    let cell_coords = SimplexCoords::from_simplex_and_coords(topology, self.mesh_coords);

    let dim = self.source.dim_intrinsic();
    let grade = self.grade();
    let dof_simps: Vec<_> = standard_subsimps(dim, grade).collect();
    let whitneys: Vec<_> = dof_simps
      .iter()
      .cloned()
      .map(|dof_simp| WhitneyLsf::standard(dim, dof_simp))
      .collect();

    let inner = multi_gramian(&geometry.to_metric_tensor().inverse(), grade);

    let mut elvec = ElVec::zeros(whitneys.len());
    for (iwhitney, whitney) in whitneys.iter().enumerate() {
      let inner_pointwise = |local: CoordRef| {
        let global = cell_coords.local2global(local);
        inner.inner(
          self
            .source
            .at_point(&global)
            .precompose_form(&cell_coords.linear_transform())
            .coeffs(),
          whitney.at_point(local).coeffs(),
        )
      };
      let value = self.qr.integrate_local(&inner_pointwise, geometry.vol());
      elvec[iwhitney] = value;
    }
    elvec
  }
}
```

## 5.3   Assembly

The element matrix provider provides the exterior grade of the bilinear form arguments. Based on this the corresponding dimension of simplices are used to assemble. The assembly process itself is a standard FEM technique [1].

We use rayon [20] to parallelize the assembly process over all the different cells, since these computations are always independent.

```rust
pub type GalMat = CooMatrix;
pub fn assemble_galmat(
  topology: &Complex,
  geometry: &MeshLengths,
  elmat: impl ElMatProvider,
) -> GalMat {
  let row_grade = elmat.row_grade();
  let col_grade = elmat.col_grade();
  let nsimps_row = topology.skeleton(row_grade).len();
  let nsimps_col = topology.skeleton(col_grade).len();

  let triplets: Vec<(usize, usize, f64)> = topology
    .cells()
```

```
    .handle_iter()
    // parallelization via rayon
    .par_bridge()
    .flat_map(|cell| {
      let geo = geometry.simplex_lengths(cell);
      let elmat = elmat.eval(&geo);

      let row_subs: Vec<_> = cell.mesh_subsimps(row_grade).collect();
      let col_subs: Vec<_> = cell.mesh_subsimps(col_grade).collect();

      let mut local_triplets = Vec::new();
      for (ilocal, &iglobal) in row_subs.iter().enumerate() {
        for (jlocal, &jglobal) in col_subs.iter().enumerate() {
          let val = elmat[(ilocal, jlocal)];
          if val != 0.0 {
            local_triplets.push((iglobal.kidx(), jglobal.kidx(), val));
          }
        }
      }

      local_triplets
    })
    .collect();
  let (rows, cols, values) = triplets.into_iter().multiunzip();
  GalMat::try_from_triplets(nsimps_row, nsimps_col, rows, cols, values).unwrap()
}
```

We also have an assembly algorithm for Galerkin vectors.

```
pub type GalVec = Vector;
pub fn assemble_galvec(
  topology: &Complex,
  geometry: &MeshLengths,
  elvec: impl ElVecProvider,
) -> GalVec {
  let grade = elvec.grade();
  let nsimps = topology.skeleton(grade).len();

  let entries: Vec<(usize, f64)> = topology
    .cells()
    .handle_iter()
    // parallelization via rayon
    .par_bridge()
    .flat_map(|cell| {
      let geo = geometry.simplex_lengths(cell);
      let elvec = elvec.eval(&geo, &cell);
      let subs: Vec<_> = cell.mesh_subsimps(grade).collect();

      let mut local_entires = Vec::new();
      for (ilocal, &iglobal) in subs.iter().enumerate() {
        if elvec[ilocal] != 0.0 {
          local_entires.push((iglobal.kidx(), elvec[ilocal]));
        }
      }
      local_entires
    })
    .collect();

  let mut galvec = Vector::zeros(nsimps);
  for (irow, val) in entries {
    galvec[irow] += val;
  }
  galvec
}
```

# Chapter 6

# Hodge-Laplacian

In this chapter we now solve some PDEs based on the Hodge-Laplace operator [4]. We consider the Hodge-Laplace eigenvalue problem and the Hodge-Laplace source problem (analog of Poisson equation).

The Hodge-Laplace operator generalizes the scalar (0-form) Laplace-Beltrami operator, to an operator acting on any differential $k$-form [5], [4]. As such the 0-form Hodge-Laplacian $\Delta^0$ is exactly the Laplace-Beltrami operator and we can write it using the exterior derivative $\mathrm{d}$ and the codifferential $\delta$.

$$\Delta^0 f = -\mathrm{div}\,\mathbf{grad}\, f = \delta^1 \mathrm{d}^0 f \tag{74}$$

The $k$-form Hodge-Laplacian $\Delta^k$ is defined as [4]

$$\begin{aligned}
\Delta^k &: \Lambda^k(\Omega) \to \Lambda^k(\Omega) \\
\Delta^k &= \mathrm{d}^{k-1}\delta^k + \delta^{k+1}\mathrm{d}^k
\end{aligned} \tag{75}$$

## 6.1 Eigenvalue Problem

We first consider the Eigenvalue problem, because it's simpler than the source problem. Furthermore the source problem *relies* on the eigenvalue problem.

The strong primal form of the Hodge-Laplace eigenvalue problem is [4]
Find $\lambda \in \mathbb{R}, u \in \Lambda^k(\Omega)$, s.t.

$$(\mathrm{d}\delta + \delta\mathrm{d})u = \lambda u \tag{76}$$

In FEM we don't solve the PDE based on the strong form, but instead we rely on a weak variational form. However the primal weak form is not suited for discretization, so instead we make use of a **mixed variational form** that includes an **auxiliary variable** $\sigma$ [3], [4].

The mixed weak form is [3], [4]
Find $\lambda \in \mathbb{R}$, $(\sigma, u) \in \left(H\Lambda^{k-1} \times H\Lambda^k \setminus \{0\}\right)$, s.t.

$$\begin{aligned}
\langle \sigma, \tau \rangle - \langle u, \mathrm{d}\tau \rangle &= 0 & \forall \tau \in H\Lambda^{k-1} \\
\langle \mathrm{d}\sigma, v \rangle + \langle \mathrm{d}u, \mathrm{d}v \rangle &= \lambda \langle u, v \rangle & \forall v \in H\Lambda^k
\end{aligned} \tag{77}$$

This formulation involves exactly the bilinear forms, we have implemented in one of the previous chapters.

We now perform Galerkin discretization of this variational problem by choosing as finite dimensional subspace of our function space $H\Lambda^k$ the space of Whitney forms $\mathcal{W}\lambda^k \subseteq H\Lambda^k$ [3] and as basis the Whitney basis $\left\{\varphi_i^k\right\}$ [10], [2]. We then replace $\sigma$ and $u$ by basis expansions $\sigma = \sum_j \sigma_j \varphi_j^{k-1}$, $u = \sum_i u_i \varphi_i^k$ and arrive at the linear system of equations.

$$\begin{aligned}
\sum_j \sigma_j \langle \varphi_j^{k-1}, \varphi_i^{k-1} \rangle - \sum_j u_j \langle \varphi_j^k, \mathrm{d}\varphi_i^{k-1} \rangle &= 0 \\
\sum_j \sigma_j \langle \mathrm{d}\varphi_j^{k-1}, \varphi_i^k \rangle + \sum_j u_j \langle \mathrm{d}\varphi_j^k, \mathrm{d}\varphi_i^k \rangle &= \lambda \sum_j u_j \langle \varphi_j^k, \varphi_i^k \rangle
\end{aligned} \tag{78}$$

We can now insert our Galerkin matrices.

$$\mathbf{M}^{k-1} = \left[\langle \varphi_i^{k-1}, \varphi_j^{k-1}\rangle\right]_{ij} \quad \mathbf{C} = \left[\langle \mathrm{d}\varphi_i^{k-1}, \varphi_j^k\rangle\right]_{ij}$$

$$\mathbf{D} = \left[\langle \varphi_i^k, \mathrm{d}\varphi_j^{k-1}\rangle\right]_{ij} \quad \mathbf{L} = \left[\langle \mathrm{d}\varphi_i^k, \mathrm{d}\varphi_j^k\rangle\right]_{ij} \quad \mathbf{M}^k = \left[\langle \varphi_i^k, \varphi_j^k\rangle\right]_{ij} \tag{79}$$

And arrive at the LSE.

$$\mathbf{M}^{k-1}\boldsymbol{\sigma} - \mathbf{C}\boldsymbol{u} = \mathbf{0}$$
$$\mathbf{D}\boldsymbol{\sigma} + \mathbf{L}\boldsymbol{u} = \lambda \mathbf{M}^k \boldsymbol{u} \tag{80}$$

This LSE has block structure and can be written as

$$\begin{bmatrix} \mathbf{M}^{k-1} & -\mathbf{C} \\ \mathbf{D} & \mathbf{L} \end{bmatrix} \begin{bmatrix} \boldsymbol{\sigma} \\ \boldsymbol{u} \end{bmatrix} = \lambda \begin{bmatrix} \mathbf{0}_{N^{k-1}\times N^{k-1}} & \mathbf{0}_{N^{k-1}\times N^k} \\ \mathbf{0}_{N^k \times N^{k-1}} & \mathbf{M}^k \end{bmatrix} \begin{bmatrix} \boldsymbol{\sigma} \\ \boldsymbol{u} \end{bmatrix} \tag{81}$$

This is a symmetric indefinite sparse generalized matrix eigenvalue problem, that can be solved by an iterative eigensolver such as Krylov-Schur [4]. In SLEPc [17] terminology this is called a GHIEP problem.

We have a helper struct for computing all the the relevant Galerkin matrices, as well as a system matrix for the mixed problem. The computation of the Galerkin matrices is done efficiently, by only assembling as few mass matrices as possible and then combining them with the full exterior derivatives, instead of assembling each bilinear form separately.

```rust
pub struct MixedGalmats {
  mass_sigma: GalMat,
  dif_sigma: GalMat,
  codif_u: GalMat,
  codifdif_u: GalMat,
  mass_u: GalMat,
}
impl MixedGalmats {
  pub fn compute(topology: &Complex, geometry: &MeshLengths, grade: ExteriorGrade) -> Self {
    let dim = topology.dim();
    assert!(grade <= dim);

    let mass_u = assemble_galmat(topology, geometry, HodgeMassElmat::new(dim, grade));
    let mass_u_csr = CsrMatrix::from(&mass_u);

    let (mass_sigma, dif_sigma, codif_u) = if grade > 0 {
      let mass_sigma = assemble_galmat(topology, geometry, HodgeMassElmat::new(dim, grade - 1));

      let exdif_sigma = topology.exterior_derivative_operator(grade - 1);
      let exdif_sigma = CsrMatrix::from(&exdif_sigma);

      let dif_sigma = &mass_u_csr * &exdif_sigma;
      let dif_sigma = CooMatrix::from(&dif_sigma);

      let codif_u = &exdif_sigma.transpose() * &mass_u_csr;
      let codif_u = CooMatrix::from(&codif_u);

      (mass_sigma, dif_sigma, codif_u)
    } else {
      (GalMat::new(0, 0), GalMat::new(0, 0), GalMat::new(0, 0))
    };

    let codifdif_u = if grade < topology.dim() {
      let mass_plus = assemble_galmat(topology, geometry, HodgeMassElmat::new(dim, grade + 1));
      let mass_plus = CsrMatrix::from(&mass_plus);
      let exdif_u = topology.exterior_derivative_operator(grade);
      let exdif_u = CsrMatrix::from(&exdif_u);
      let codifdif_u = exdif_u.transpose() * mass_plus * exdif_u;
      CooMatrix::from(&codifdif_u)
    } else {
      GalMat::new(0, 0)
    };
```

```rust
  Self {
    mass_sigma,
    dif_sigma,
    codif_u,
    codifdif_u,
    mass_u,
  }
}
pub fn sigma_len(&self) -> usize { self.mass_sigma.nrows() }
pub fn u_len(&self) -> usize { self.mass_u.nrows() }

pub fn mixed_hodge_laplacian(&self) -> CooMatrix {
  let Self {
    mass_sigma,
    dif_sigma,
    codif_u,
    codifdif_u,
    ..
  } = self;
  let codif_u = codif_u.clone();
  CooMatrix::block(&[&[mass_sigma, &(codif_u.neg())], &[dif_sigma, codifdif_u]])
}
}
```

The code for solving the actual EVP is then very simple. It calls the PETSc/SLEPc solver.

```rust
pub fn solve_hodge_laplace_evp(
  topology: &Complex,
  geometry: &MeshLengths,
  grade: ExteriorGrade,
  neigen_values: usize,
) -> (Vector, Matrix, Matrix) {
  let galmats = MixedGalmats::compute(topology, geometry, grade);

  let lhs = galmats.mixed_hodge_laplacian();

  let sigma_len = galmats.sigma_len();
  let u_len = galmats.u_len();
  let mut rhs = CooMatrix::zeros(sigma_len + u_len, sigma_len + u_len);
  for (mut r, mut c, &v) in galmats.mass_u.triplet_iter() {
    r += sigma_len;
    c += sigma_len;
    rhs.push(r, c, v);
  }

  let (eigenvals, eigenvectors) = petsc_ghiep(&(&lhs).into(), &(&rhs).into(), neigen_values);

  let eigen_sigmas = eigenvectors.rows(0, sigma_len).into_owned();
  let eigen_us = eigenvectors.rows(sigma_len, u_len).into_owned();
  (eigenvals, eigen_sigmas, eigen_us)
}
```

## 6.2  Source Problem

The Hodge-Laplace Source Problem is the generalization of the Poisson equation to arbitrary differential $k$-forms [4]. In strong form it is:
Find $u \in \Lambda^k(\Omega)$, given $f \in \Lambda^k(\Omega)$, s.t.

$$\Delta u = f - P_{\mathfrak{H}} f \quad \text{in } \Omega$$
$$\text{tr} \star u = 0, \quad \text{tr}\,\mathrm{d}u = 0 \quad \text{on } \partial\Omega \tag{82}$$
$$u \perp \mathfrak{H}$$

This equation is not quite as simple as the normal Poisson equation $\Delta u = f$. Instead it includes two additional parts involving $\mathfrak{H}$, which is the space of harmonic forms $\mathfrak{H}^k = \ker \Delta = \{v \in \Lambda^k \mid \Delta v = 0\}$ [4]. The first change is that we remove the harmonic part $P_{\mathfrak{H}} f$ of $f$. The second difference is that we require that our solution $u$ is orthogonal to harmonic forms.

Harmonic forms are concrete representative of the cohomology quotient group $\mathcal{H}^k = \frac{\ker \mathrm{d}}{\mathrm{im}\, \mathrm{d}}$ [4], [5].

We once again tackle a mixed weak formulation based on the auxiliary variable $\sigma$ and this time a second one $p$ that represents $f$ without harmonic component [3], [4].

Given $f \in L^2\Lambda^k$, find $(\sigma, u, p) \in \left(H\Lambda^{k-1} \times H\Lambda^k \times \mathfrak{H}^k\right)$ s.t.

$$
\begin{aligned}
\langle \sigma, \tau \rangle - \langle u, \mathrm{d}\tau \rangle &= 0 & \forall \tau \in H\Lambda^{k-1} \\
\langle \mathrm{d}\sigma, v \rangle + \langle \mathrm{d}u, \mathrm{d}v \rangle + \langle p, v \rangle &= \langle f, v \rangle & \forall v \in H\Lambda^k \\
\langle u, q \rangle &= 0 & \forall q \in \mathfrak{H}^k
\end{aligned}
\tag{83}
$$

We once again perform Galerkin discretization [3].

$$
\begin{aligned}
\sum_j \sigma_j \langle \varphi_j^{k-1}, \varphi_i^{k-1} \rangle - \sum_j u_j \langle \varphi_j^k, \mathrm{d}\varphi_i^{k-1} \rangle &= 0 \\
\sum_j \sigma_j \langle \mathrm{d}\varphi_j^{k-1}, \varphi_i^k \rangle + \sum_j u_j \langle \mathrm{d}\varphi_j^k, \mathrm{d}\varphi_i^k \rangle + \sum_j p_j \langle \eta_j^k, \varphi_i^k \rangle &= \langle f, \varphi_i^k \rangle \\
\sum_j u_j \langle \varphi_j^k, \eta_i^k \rangle &= 0
\end{aligned}
\tag{84}
$$

By inserting our known Galerkin matrices, we obtain.

$$
\begin{aligned}
\mathbf{M}^{k-1}\boldsymbol{\sigma} - \mathbf{C}\boldsymbol{u} &= 0 \\
\mathbf{D}\boldsymbol{\sigma} + \mathbf{L}\boldsymbol{u} + \mathbf{MH}\boldsymbol{p} &= \boldsymbol{b} \\
\mathbf{H}^\top \mathbf{M}\boldsymbol{u} &= 0
\end{aligned}
\tag{85}
$$

Or in block-structure

$$
\begin{bmatrix}
\mathbf{M}^{k-1} & -\mathbf{C} & 0 \\
\mathbf{D} & \mathbf{L} & \mathbf{MH} \\
0 & \mathbf{H}^\top \mathbf{M} & 0
\end{bmatrix}
\begin{bmatrix}
\boldsymbol{\sigma} \\
\boldsymbol{u} \\
\boldsymbol{p}
\end{bmatrix}
=
\begin{bmatrix}
0 \\
\boldsymbol{b} \\
0
\end{bmatrix}
\tag{86}
$$

Where the right-hand side $\boldsymbol{b}$, corresponding to the source term, is approximated via quadrature, as previously discussed.

Computing harmonics is just a matter of solving the EVP to obtain the eigenfunctions that correspond to eigenvalue zero.

```rust
pub fn solve_hodge_laplace_harmonics(
  topology: &Complex,
  geometry: &MeshLengths,
  grade: ExteriorGrade,
  homology_dim: usize,
) -> Matrix {
  if homology_dim == 0 {
    let nwhitneys = topology.nsimplices(grade);
    return Matrix::zeros(nwhitneys, 0);
  }

  let (eigenvals, _, harmonics) = solve_hodge_laplace_evp(topology, geometry, grade, homology_dim);
  assert!(eigenvals.iter().all(|&eigenval| eigenval <= 1e-12));
  harmonics
}
```

To solve the actual source problem, we now just need to assemble the system matrix and call the PETSc solver.

```rust
pub fn solve_hodge_laplace_source(
  topology: &Complex,
  geometry: &MeshLengths,
```

```
  source_galvec: GalVec,
  grade: ExteriorGrade,
  homology_dim: usize,
) -> (Cochain, Cochain, Cochain) {
  let harmonics = solve_hodge_laplace_harmonics(topology, geometry, grade, homology_dim);

  let galmats = MixedGalmats::compute(topology, geometry, grade);

  let mass_u = CsrMatrix::from(&galmats.mass_u);
  let mass_harmonics = &mass_u * &harmonics;

  let sigma_len = galmats.sigma_len();
  let u_len = galmats.u_len();

  let mut galmat = galmats.mixed_hodge_laplacian();

  galmat.grow(mass_harmonics.ncols(), mass_harmonics.ncols());

  for (mut r, mut c) in (0..mass_harmonics.nrows()).cartesian_product(0..mass_harmonics.ncols()) {
    let v = mass_harmonics[(r, c)];
    r += sigma_len;
    c += sigma_len + u_len;
    galmat.push(r, c, v);
  }
  for (mut r, mut c) in (0..mass_harmonics.nrows()).cartesian_product(0..mass_harmonics.ncols()) {
    let v = mass_harmonics[(r, c)];
    // transpose
    mem::swap(&mut r, &mut c);
    r += sigma_len + u_len;
    c += sigma_len;
    galmat.push(r, c, v);
  }

  let system_matrix = CsrMatrix::from(&galmat);

  #[allow(clippy::toplevel_ref_arg)]
  let rhs = na::stack![
    Vector::zeros(sigma_len);
    source_galvec;
    Vector::zeros(harmonics.ncols());
  ];

  let galsol = petsc_saddle_point(&system_matrix, &rhs);
  let sigma = Cochain::new(grade - 1, galsol.view_range(..sigma_len, 0).into_owned());
  let u = Cochain::new(
    grade,
    galsol
      .view_range(sigma_len..sigma_len + u_len, 0)
      .into_owned(),
  );
  let p = Cochain::new(
    grade,
    galsol.view_range(sigma_len + u_len.., 0).into_owned(),
  );
  (sigma, u, p)
}
```

# Chapter 7

# Results

In this chapter, we present numerical results to verify the functionality and validate the implementation of the `formoniq` library.

Specifically, we examine the Hodge-Laplacian eigenvalue problem on a domain with non-trivial topology and perform a convergence study for the Hodge-Laplacian source problem using the Method of Manufactured Solutions on a $n$-dimensional hypercube domain.

## 7.1  1-form Eigenvalue Problem on Torus

To assess the library's ability to handle **non-trivial topologies** as well as **globally curved geometry**, we solve a eigenvalue problem to compute the spectrum of the 1-form Hodge-Laplacian $\Delta^1$ on a torus $\mathbb{T}^2$.

A mesh of a torus with major radius $r_1 = 0.5$ and minor radius $r_2 = 0.2$ was generated using Gmsh [21], by specifying a `.geo` file with the line:

```
Torus(1) = {-0, -0, 0, 0.5, 0.2, 2*Pi};
```

The topology of the torus is characterized by its first **Betti number** $b_1 = 2$, which predicts a two-dimensional kernel for $\Delta^1$ spanned by harmonic 1-forms representing the two fundamental cycles of the domain.

The theoretical eigenvalues of $\Delta^1$ on an idealized flat torus with circumferences $L_1 = 2\pi r_1 = \pi$ and $L_2 = 2\pi r_2 = 0.4\pi$ are given by $\lambda_{m,n} = \left(\frac{2\pi m}{L_1}\right)^2 + \left(\frac{2\pi n}{L_2}\right)^2 = 4m^2 + 25n^2$ for integers $m, n \in \mathbb{Z}$.

The lowest computed eigenvalues obtained using formoniq and the SLEPc [17] eigensolver are:

```
ieigen=0, eigenval=-0.000
ieigen=1, eigenval=0.000
ieigen=2, eigenval=4.116
ieigen=3, eigenval=4.116
ieigen=4, eigenval=4.116
ieigen=5, eigenval=4.116
ieigen=6, eigenval=14.447
ieigen=7, eigenval=14.447
ieigen=8, eigenval=14.449
ieigen=9, eigenval=14.449
ieigen=10, eigenval=24.649
```

These numerical results show excellent agreement with theoretical expectations:

- The computed spectrum correctly identifies the two zero eigenvalues $\lambda = 0.000$ corresponding to the two harmonic 1-forms $m = n = 0$, accurately capturing the torus's topology ($b_1 = 2$).
- The first non-zero eigenvalue group is computed as $\lambda \approx 4.116$ with multiplicity 4. This closely matches the theoretical value $\lambda_{1,0} = 4$ and its expected multiplicity.
- The second non-zero group is computed around $\lambda \approx 14.447 \approx 14.449$ with multiplicity 4, corresponding well to the theoretical value $\lambda_{2,0} = 16$ and its expected multiplicity.
- The next computed eigenvalue $\lambda \approx 24.649$ aligns closely with the theoretical value $\lambda_{0,1} = 25$.

Figure 4 shows visualizations of the two harmonic 1-forms, represented by their vector proxies.

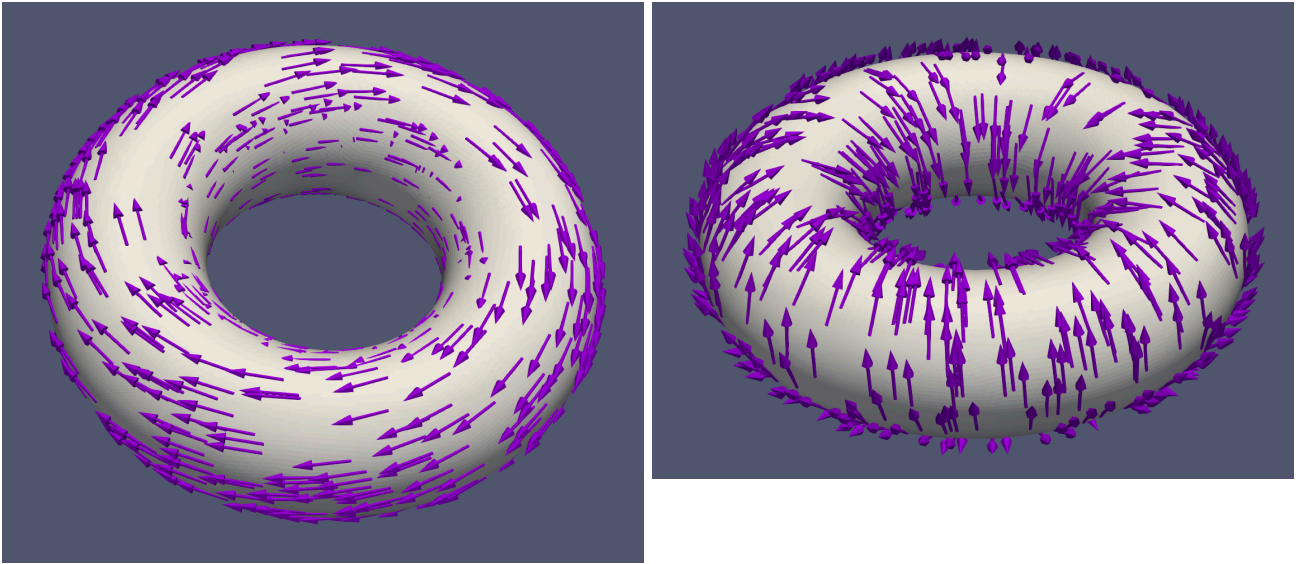

Figure 4: The two harmonic forms on the torus, representing the 1-cohomology spaces.

The slight deviations between computed values and the ideal flat torus eigenvalues are expected due to discretization error. However, the accurate recovery of the zero eigenvalues and the correct multiplicities for the lowest eigenvalue groups strongly validates the FEEC implementation for eigenvalue problems on domains with non-trivial topology.

This example program which computes the eigenvalue works on any mesh. The program asks for the `.obj` or `.msh` (gmsh) file, then loads it in and computes the eigenpairs. The FE eigenfunction is then sampled at each barycenter and the values are written into a file, which can then be visualized, with for instance paraview.

The example program can be run using

```
cargo run --release --example hodge_laplace_source
```

## 7.2  Source Problem

We verify the source problem by means of the **method of manufactured solution** [1]. We restrict ourselves to a simple setup in globally flat geometry on a subset of $\mathbb{R}^n$. We validate only the 1-form source problem, but we do this in arbitrary dimensions $n \geq 2$.

Our domain is

$$\Omega = [0, \pi]^n \tag{87}$$

For our manufactured solution, we've chosen a simple 1-form $u \in \Lambda^1(\Omega)$ that generalizes easily to arbitrary dimensions.

$$u = \sum_{i=1}^{n} u_i \mathrm{d}x^i \quad \text{with} \quad u_i = \sin^2(x^i) \prod_{j \neq i} \cos(x^j) \tag{88}$$

When using vector proxies $u^\sharp$ we get the following vector fields in the 2D and 3D case.

$$n = 2 \implies u^\sharp = \begin{bmatrix} \sin^2(x)\cos(y) \\ \cos(x)\sin^2(y) \end{bmatrix} \qquad n = 3 \implies u^\sharp = \begin{bmatrix} \sin^2(x)\cos(y)\cos(z) \\ \cos(x)\sin^2(y)\cos(z) \\ \cos(x)\cos(y)\sin^2(z) \end{bmatrix} \tag{89}$$

We can analytically derive the exterior derivative $\mathrm{d}u$ of this solution $u$.

$$\mathrm{d}u = \sum_{k<i} \left[ \left( \prod_{j \neq i,k} \cos(x^j) \right) \sin(x^i)\sin(x^k)\big(\sin(x^k) - \sin(x^i)\big) \right] \mathrm{d}x^k \wedge \mathrm{d}x^i \tag{90}$$

In the method of manufactured solution, the source term $f$ is set to equal the 1-form Hodge-Laplacian $\Delta^1 u$ of the exact solution $u$.

$$f = \Delta^1 u \tag{91}$$

The corresponding source term $f = \Delta^1 u$ is computed analytically:

$$(\Delta^1 \boldsymbol{u})_i = \Delta^0 u_i = -\left(2\cos(2x^i) - (n-1)\sin^2(x^i)\right) \prod_{j \neq i} \cos(x^j) \tag{92}$$

Our solution and it's derivative have boundary traces that are equal to zero. This leads to homogeneous natural boundary conditions, meaning no additional terms are required in the variational formulation.

$$\text{Tr}_{\partial\Omega} u = 0 \qquad \text{Tr}_{\partial\Omega} \mathrm{d}u = 0 \tag{93}$$

We use formoniq to solve this problem and determine the rate of convergence for dimensions 2 and 3, but even higher dimensions would work.

For each dimension we generate finer and finer meshes with mesh width $h$ halved in each step. The meshes are generated using our custom tensor-product triangulation algorithm. We compute the right hand side vector by assembling the element matrix provider for the source term based on the exact Laplacian. Then we just call the `solve_hodge_laplace_source` routine.

We compute the $L^2$ norm of the error in the function value $\|u - u_h\|_{L^2}$. We do this by computing the pointwise error norm $\|u - u_h\|_{\Lambda^k}$ and integrating it using quadrature of order 3, which is sufficient for the quadrature error to not dominate the finite element error [1]. We also compute the $L^2$ error in the exterior derivative, $\|\mathrm{d}u - \mathrm{d}u_h\|_{L^2}$ using the same approach.

```rust
pub fn fe_l2_error<E: ExteriorField>(
  fe_cochain: &Cochain,
  exact: &E,
  topology: &Complex,
  coords: &MeshCoords,
) -> f64 {
  let dim = topology.dim();
  let qr = SimplexQuadRule::order3(dim);
  let fe_whitney = WhitneyForm::new(fe_cochain.clone(), topology, coords);
  let inner = multi_gramian(&Gramian::standard(dim), fe_cochain.dim());
  let error_pointwise = |x: CoordRef, cell: SimplexHandle| {
    inner.norm_sq((exact.at_point(x) - fe_whitney.eval_known_cell(cell, x)).coeffs())
  };
  qr.integrate_mesh(&error_pointwise, topology, coords).sqrt()
}
```

This is the code for our convergence test. It can be run using `cargo run --release --example hodge_laplace_source`.

```rust
fn main() -> Result<(), Box<dyn std::error::Error>> {
  tracing_subscriber::fmt::init();
  let path = "out/laplacian_source";
  let _ = fs::remove_dir_all(path);
  fs::create_dir_all(path).unwrap();

  let grade = 1;
  let homology_dim = 0;

  for dim in 2_usize..=3 {
    println!("Solving Hodge-Laplace in {dim}d.");

    let solution_exact = DiffFormClosure::one_form(
      |p| {
        Vector::from_iterator(
          p.len(),
          (0..p.len()).map(|i| {
            let prod = p.remove_row(i).map(|a| a.cos()).product();
            p[i].sin().powi(2) * prod
          }),
        )
      },
      dim,
```

```rust
);

let dif_solution_exact = DiffFormClosure::new(
  Box::new(move |p: CoordRef| {
    let dim = p.len();
    let ncomponents = if dim > 1 { dim * (dim - 1) / 2 } else { 0 };
    let mut components = Vec::with_capacity(ncomponents);

    let sin_p: Vec<_> = p.iter().map(|&pi| pi.sin()).collect();
    let cos_p: Vec<_> = p.iter().map(|&pi| pi.cos()).collect();

    for k in 0..dim {
      for i in (k + 1)..dim {
        let mut prod_cos = 1.0;
        #[allow(clippy::needless_range_loop)]
        for j in 0..dim {
          if j != i && j != k {
            prod_cos *= cos_p[j];
          }
        }
        let coeff = prod_cos * sin_p[i] * sin_p[k] * (sin_p[k] - sin_p[i]);
        components.push(coeff);
      }
    }
    ExteriorElement::new(components.into(), dim, 2)
  }),
  dim,
  2,
);

let laplacian_exact = DiffFormClosure::one_form(
  |p| {
    Vector::from_iterator(
      p.len(),
      (0..p.len()).map(|i| {
        let prod: f64 = p.remove_row(i).map(|a| a.cos()).product();
        -(2.0 * (2.0 * p[i]).cos() - (p.len() - 1) as f64 * p[i].sin().powi(2)) * prod
      }),
    )
  },
  dim,
);

println!(
  "| {:>2} | {:8} | {:>7} | {:>8} | {:>7} |",
  "k", "L2 err", "L2 conv", "Hd err", "Hd conv",
);

let mut errors_l2 = Vec::new();
let mut errors_hd = Vec::new();
for irefine in 0..=(15 / dim as u32) {
  let refine_path = &format!("{path}/refine{irefine}");
  fs::create_dir_all(refine_path).unwrap();

  let nboxes_per_dim = 2usize.pow(irefine);
  let box_mesh = CartesianMeshInfo::new_unit_scaled(dim, nboxes_per_dim, PI);
  let (topology, coords) = box_mesh.compute_coord_complex();
  let metric = coords.to_edge_lengths(&topology);

  let source_data = assemble_galvec(
    &topology,
    &metric,
    SourceElVec::new(&laplacian_exact, &coords, None),
  );
```

```rust
    let (_, galsol, _) = hodge_laplace::solve_hodge_laplace_source(
      &topology,
      &metric,
      source_data,
      grade,
      homology_dim,
    );

    let conv_rate = |errors: &[f64], curr: f64| {
      errors
        .last()
        .map(|&prev| algebraic_convergence_rate(curr, prev))
        .unwrap_or(f64::INFINITY)
    };

    let error_l2 = fe_l2_error(&galsol, &solution_exact, &topology, &coords);
    let conv_rate_l2 = conv_rate(&errors_l2, error_l2);
    errors_l2.push(error_l2);

    let dif_galsol = galsol.dif(&topology);
    let error_hd = fe_l2_error(&dif_galsol, &dif_solution_exact, &topology, &coords);
    let conv_rate_hd = conv_rate(&errors_hd, error_hd);
    errors_hd.push(error_hd);

    println!(
      "| {:>2} | {:<8.2e} | {:>7.2} | {:<8.2e} | {:>7.2} |",
      irefine, error_l2, conv_rate_l2, error_hd, conv_rate_hd
    );
  }
}

  Ok(())
}
```

The output is

```
Solving Hodge-Laplace in 2d.
| k | L2 err  | L2 conv | Hd err  | Hd conv |
| 0 | 1.96e0  |     inf | 6.51e-1 |     inf |
| 1 | 1.57e0  |    0.31 | 4.31e-1 |    0.60 |
| 2 | 8.02e-1 |    0.97 | 3.51e-1 |    0.30 |
| 3 | 4.03e-1 |    0.99 | 1.35e-1 |    1.37 |
| 4 | 1.99e-1 |    1.02 | 6.00e-2 |    1.17 |
| 5 | 9.91e-2 |    1.01 | 2.89e-2 |    1.05 |
| 6 | 4.95e-2 |    1.00 | 1.43e-2 |    1.01 |
| 7 | 2.48e-2 |    1.00 | 7.15e-3 |    1.00 |
Solving Hodge-Laplace in 3d.
| k | L2 err  | L2 conv | Hd err  | Hd conv |
| 0 | 3.66e0  |     inf | 1.09e0  |     inf |
| 1 | 2.56e0  |    0.52 | 1.76e0  |   -0.69 |
| 2 | 1.46e0  |    0.80 | 7.49e-1 |    1.23 |
| 3 | 7.71e-1 |    0.93 | 3.08e-1 |    1.28 |
| 4 | 3.85e-1 |    1.00 | 1.39e-1 |    1.15 |
| 5 | 1.92e-1 |    1.00 | 6.73e-2 |    1.04 |
```

We observe experimental convergence $\mathcal{O}(h^1)$ for the $L^2$ error of the exterior derivative $\|\mathrm{d}u - \mathrm{d}u_h\|_{L^2}$ and converge $\mathcal{O}(h^1)$ for the $L^2$ error of the solution value itself $\|u - u_h\|_{L^2}$.

The observed rate $\alpha = 1$ for the exterior derivative is promising. It matches the convergence rate expected for the $\|\mathrm{d} \cdot\|_{L^2}$ component within the $H\Lambda$ norm for first-order Whitney elements, as predicted by theory [4]. This suggests a valid implementation of the discretization related to the exterior derivative operator.

The interpretation of the $L^2$ rate $\alpha = 1$ for the function value is less straightforward. While this rate is compatible with an overall $\mathcal{O}(h^1)$ convergence in the full $H\Lambda$ norm, standard FEEC theory often predicts $\mathcal{O}(h^2)$ convergence for the $L^2$ error itself via an Aubin-Nitsche duality argument [4]. The reason why this higher rate is not observed in our results requires further investigation, potentially related to the quadrature rules used for error computation or other implementation details.

It should be noted that this study constitutes only a partial validation of the implementation. A complete verification would require measuring the error in the codifferential, $\|\delta(u - u_h)\|_{L^2}$, to assess convergence in the the full $H\Lambda$ norm or the associated energy norm.

$$\|\cdot\|_{H\Lambda} = \|\cdot\|_{L^2} + \|\cdot\|_{H(\mathrm{d})} + \|\cdot\|_{H(\delta)} \tag{94}$$

This analysis was not performed here, due to time constraints.

In Figure 5 we provide a visualization of the 2D finite element solution at refinement level 4 our library has produced in the form of a vector field proxy together with a heat map of the magnitude.

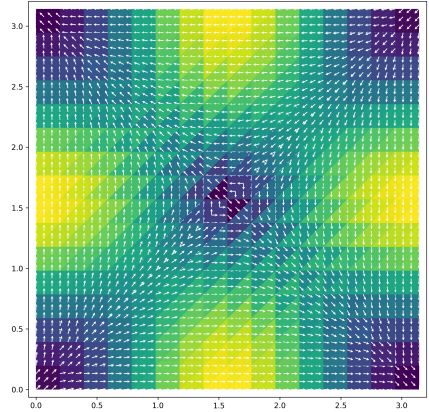

Figure 5: Finite element solution to manufactured source problem in 2D at refinement level 4.

# Chapter 8

# Conclusion and Outlook

This thesis presented `formoniq`, a novel implementation of Finite Element Exterior Calculus (FEEC) developed in the Rust programming language. The core contribution lies in its foundation on intrinsic, coordinate-free geometry and its capability to operate on simplicial complexes in arbitrary dimensions. By leveraging Rust's performance and safety features, we aimed to provide a modern and robust tool for structure-preserving discretization of partial differential equations formulated in the language of differential forms.

We successfully developed modules for handling the topology of arbitrary-dimensional simplicial complexes, representing intrinsic Riemannian geometry via edge lengths inspired by Regge Calculus, performing exterior algebra computations, and implementing discrete differential forms using cochains and first-order Whitney basis functions. Building upon this foundation, we implemented the necessary Galerkin operators for FEEC, specifically targeting the mixed weak formulation of the Hodge-Laplace equation.

`formoniq` demonstrates the feasibility and benefits of combining the rigorous mathematical framework of intrinsic FEEC with the modern software engineering practices enabled by Rust.

## 8.1 Future Work

Several avenues exist for extending and enhancing the `formoniq` library:

- Higher-Order FEEC: Extend the implementation to support higher-order polynomial basis functions for differential forms. This would enable higher accuracy and faster convergence for problems with smooth solutions but must be coupled with corresponding higher-order representations of the manifold geometry (curved simplices) to avoid introducing non-admissible geometric variational crimes. [4], [2], [11]
- Maxwell's Equations: Apply the framework to solve Maxwell's equations, particularly in contexts where the coordinate-free and arbitrary-dimensional nature is advantageous, such as electromagnetism on curved spacetimes. FEEC provides a natural and structure-preserving discretization for these equations. [26], [3], [5]
- Optimization: Profile and optimize core computational routines in Rust, potentially exploring alternative sparse matrix libraries or parallelization strategies beyond the assembly loop currently handled by Rayon.

## 8.2 Comparison to Other Implementations

The field of computational differential geometry and structure-preserving discretizations has seen several software development efforts. To position `formoniq` within this landscape, we briefly compare it to some notable existing libraries based on available documentation and publications. `formoniq` did not directly draw inspiration from these specific implementations but addresses similar challenges.

**Formoniq (This Thesis)**
- **Focus:** FEEC
- **Dimension:** Arbitrary
- **Mesh:** Simplicial Complexes
- **Geometry:** Intrinsic Regge metric
- **Discretization:** 1st order Whitney forms
- **Language:** Rust
- **Key Features:** Emphasis on FEEC on coordinate-free intrinsic geometry in arbitrary dimensions.

- **Repository:** github:luiswirth/formoniq

**PyDEC**
- **Focus:** Primarily DEC, some FEEC elements
- **Dimension:** Arbitrary
- **Mesh:** Simplicial and Cubical Complexes
- **Geometry:** Embedded
- **Discretization:** Cochains
- **Language:** Python
- **Key Features:** Mature library for DEC. Includes tools for cohomology and Hodge decomposition. [27]
- **Repository:** github:hirani/pydec

**FEEC++ / simplefem**
- **Focus:** FEEC
- **Dimension:** Hardcoded 1D, 2D and 3D
- **Mesh:** Simplicial Complexes
- **Geometry:** Embedded
- **Discretization:** Arbitrary order polynomial differential forms
- **Language:** C++
- **Key Features:** Focus on arbitrary polynomial order differential forms, including Whitney and Sullivan forms. Comes with all necessary linear algebra subroutines. [28]
- **Repository:** github:martinlicht/simplefem

**DDF.jl**
- **Focus:** Foundational tools for DEC and FEEC
- **Dimension:** Arbitrary
- **Mesh:** Simplicial Complexes
- **Geometry:** Embedded
- **Discretization:** Higher-order discretizations
- **Language:** Julia
- **Key Features:** Arbitrary dimensions and higher-order methods. Unfinished. [29]
- **Repository:** github:eschnett/DDF.jl

**dexterior**
- **Focus:** DEC
- **Dimension:** Arbitrary
- **Mesh:** Simplicial Complexes
- **Geometry:** Embedded
- **Discretization:** Cochain
- **Language:** Rust
- **Key Features:** DEC in Rust, inspired by PyDEC. wgpu visualizer. [30]
- **Repository:** github:m0lentum/dexterior

This comparison highlights `formoniq`'s specific niche: providing a arbitrary-dimensional FEEC implementation fundamentally based on intrinsic geometry, currently focused on first-order methods. It complements existing libraries by offering a different language choice (Rust) and a distinct focus on the coordinate-free geometric perspective inherent in FEEC.

# Appendix A

# Formoniq Source Code

All source code of our implementation **Formoniq** is available on GitHub.

## github:luiswirth/formoniq

The main component of the Git repository is the Rust code.

In order to run the Rust code, the Rust toolchain needs to be installed on the machine. This is best done through the official Rust installer and toolchain manager Rustup, which is easily installed by running a single `curl` command in the shell, which can be found on `rustup.rs`.

In the repository root `./` we have a `./Cargo.toml` file that specifies a Cargo workspace. Therefore **all Cargo commands should be run from this root directory**. In `./crates/` all the various libraries that make up Formoniq are placed. Each of them follows the standard structure of a Rust library. In `./crates/<crate>/src/` the source code of the current library can be found.

Formoniq is first-and-foremost a library, so just running it doesn't make any sense, since it is not an executable. However we do provide various example codes that do produce executables and directly make use of the functionality of our library. All of these examples can be found under `./crates/formoniq/examples/`. These examples can be run using the following command.

## `cargo` `run --example` `<example>`

This command compiles all crates and runs the example in the shell. It is recommended to additionally add the `--release` flag to build in release mode, instead of debug mode to profit from optimizations. To figure out the names of the available examples one can look into the example directory or use tab-completion in the cargo command.

A documentation website can be generated using rustdoc, by running the following command

`cargo doc --no-deps --open`

This can serve as a index of all the datastructures we have developed.

## A.1  PETSc/SLEPc Solver

Due to the immature Rust sparse linear algebra ecosystem, we have to rely on C/C++ implementations of sparse solvers. For this we rely on **PETSc** for direct LSE solvers and on **SLEPc**, which builds on PETSc, for eigensolvers. We've written small PETSc and SLEPc programs that can be found under `./petsc-solver`. These programs load the system matrices from disk, solve the problem and write the solution back to disk. Our Rust program interoperates with these programs, by doing this writing and reading to disk as well, to communicate the problem setup and retrieve the solution. It then also automatically calls the PETSc/SLEPc program. The user does never directly interact with these solvers, but the Rust code manages this itself.

The user needs to compile this small solver program. For this both PETSc and SLEPc need to be installed on the system. We refer to the official PETSc and SLEPc documentation on [petsc.org](petsc.org) and [slepc.upv.es](slepc.upv.es), where installation of both software suites is explained. The steps that worked for the author are outlined in Section C.

Given that PETSc and SLEPc are successfully installed, we can now build our small PETSc/SLEPc solver. For this it is crucial that both `PETSC_DIR` and `SLEPC_DIR` are set.

```
export PETSC_DIR=<PATH_TO_PETSC>
export SLEPC_DIR=<PATH_TO_SLEPC>
```

Then we can navigate into the `./petsc-solver` directory and simply run make. This will produce two executables `ghiep.o` and `hils.o`. These executables will then be found and run by the Formoniq.

## A.2   Plotting

We create a simple python matplotlib script that is capable of visualizing Whitney 1-forms that correspond to arbitrary 1-cochains on arbitrary 2D meshes embedded in 2D. The necessary code can be found in `./plot` and is managed using the uv python package manager, which can be found under [github:astral-sh/uv](github:astral-sh/uv).

The usage should be clear from this `--help` message.

```
> uv run src/main.py --help
usage: main.py [-h] [--skip-zero] [--heatmap-res HEATMAP_RES]
               [--quiver-count QUIVER_COUNT] [--highlight]
               path

Whitney Forms Visualizer

positional arguments:
  path                  Path to the input files

options:
  -h, --help            show this help message and exit
  --skip-zero           Skip triangles with all zero DOFs
  --heatmap-res HEATMAP_RES
                        Resolution of the heatmap (default: 30)
  --quiver-count QUIVER_COUNT
                        Number of quiver arrows per triangle dimension
                        (default: 20)
  --highlight           Disable highlighting of non-zero DOF edges
```

In the directory `./plot/in/` various simple input are prepared for visualizing local shape functions on the reference triangle and global shape functions on the equilateral "triforce" mesh. One can for example run the following command for some nice visuals that are also shows in the thesis.

```
uv run src/main.py --quiver-count=5 --heatmap-res=10 in/triforce/
```

# Appendix B

# Thesis Typst Source Code

The thesis document itself has been written using the new type-setting language Typst, which can be found under typst.app. It is very similar to LaTeX, but aims to elevate various points of frustration that are common in LaTeX and to provide a more modern feel.

The source code for this Typst document, is available on GitHub.

<div align="center">

## github:luiswirth/bsc-thesis

</div>

To build the document in the form of a PDF, the Typst compiler needs to be installed and then the `./build.sh` script can be run from the root project directory.

# Appendix C

# PETSc/SLEPc Installation

We quickly explain the installation steps for PETSc and SLEPc that worked for the author.

## C.1 PETSc

For PETSc we simply cloned the repository, configured the installation and the ran make to compile all files.

```
git clone -b release https://gitlab.com/petsc/petsc.git petsc
cd petsc

./configure \
  --with-cc=gcc --with-cxx=g++ --with-fc=gfortran \
  --download-mpich --download-fblaslapack \
  --download-mumps --download-scalapack --download-parmetis \
  --download-metis --download-ptscotch

make all check
```

Alternatively one could have build with MPI support, but for us distributed computations did not work.

```
--with-cc=mpicc --with-cxx=mpicxx --with-fc=mpif90
```

Now PETSc should be successfully installed.

## C.2 SLEPc

SLEPc depends on PETSc and therefore needs to be installed after it and told it's location. For this we set the `PETSC_DIR` environment variable to point to the PETSc installation directory. If you are inside the PETSc directory, then the following works.

```
export PETSC_DIR=$(pwd)
```

Also `PETSC_ARCH` has to be set to the architecture for which PETSc has been build. This is the name of one of the folder that was produced in the PETSc directory. For our Linux debug build, we need to set the envvar to the following.

```
export PETSC_ARCH=arch-linux-c-debug
```

To do the actual SLEPc installation the steps are very similar to PETSc: We clone the repository, where we make sure we have the same version of SLEPc as for PETSc, then run the configure script and finally run make.

```
git clone https://gitlab.com/slepc/slepc
git checkout v3.23.0
cd slepc

./configure

make
make test
```

If everything went well without any errors, SLEPc should be successfully installed.